%% file: SWE-AVS-FE.tex
\documentclass[review]{elsarticle}

 \pdfoutput=1

\usepackage{theorem}
\usepackage{amssymb,amsmath}

\usepackage{mathptmx}       
\usepackage{helvet}         
\usepackage{courier}        
\usepackage{type1cm}        
\usepackage{bm}                            
\usepackage{makeidx}         
\usepackage{graphicx}        
\usepackage{multicol}        
\usepackage[bottom]{footmisc}
\usepackage{subfigure}
\usepackage{booktabs}

\usepackage{epsfig} 
\usepackage{epstopdf}
\usepackage{color}
\usepackage{colordvi}

\newcommand{\norm}[2]         { \| {#1} \|_{#2} }                      
\newcommand{\bfm}[1]             { \mathbf{#1}     }             %
\newcommand{\qq}             { \bfm{q}     }             
\newcommand{\Nabla}       { \boldsymbol{\nabla} }   
\newcommand{\SLTO}            { L^2(\Omega) }                       
\newcommand{\SLTOT}            { L^2(\Omega_{\text{T}}) }                       
\newcommand{\SHOO}            { H^1(\Omega) }                         
\newcommand{\SHOP}            { H^1(\Ph) }                              
\newcommand{\SHdivP}            { H(\text{div},\Ph) }            
\newcommand{\VVK}            { {V(K_m)} }              
\newcommand{\VV}            { {V(\Ph)} }        
\newcommand{\UUU}            { U(\Omega) }           
\newcommand{\UUUT}            { U(\Omega_{\text{T}}) }           
\newcommand{\UUUhT}            { U^h(\Omega_{\text{T}}) }           
\newcommand{\VVh}            { V^*(\Ph) }          
\newcommand{\SHOK}            { H^1(K_m) }                        
\newcommand{\SHdivO}            { H(\text{\textbf{div}},\Omega) }             
\newcommand{\SHdivK}            { H(\text{\textbf{div}},K_m) }            
\newcommand{\SHdivPh}            { H(\text{\textbf{div}},\Ph) }
\newcommand{\SHOOT}            { H^1(\Omega_{\text{T}}) }                   

\newcommand{\bb}             { {\bfm{b}}   }             
\newcommand{\xx}             { \bfm{x}     }             
\newcommand{\vn}             { \bfm{n}     }             

\newcommand{\Ph}            { \mathcal{P}_h }
\newcommand{\Kep}           { K_m \in \Ph}
\newcommand{\dKm}           { \partial K_m }

\newcommand{\dx}              { \; {\rm d} \bfm{x}   }                
\newcommand{\dss}             { \, {\rm d} s   }                          
\newcommand{\summa}[2]        { \overset{#2}{\underset{#1}{\sum}} } 
\newcommand{\supp}[1]         { \underset{#1}{\sup} \, }        

\newcommand{\uu}              { \bfm{u} }                      
\newcommand{\vv}              { \bfm{v} }                      
\newcommand{\www}             { \bfm{w} }                      
\newcommand{\nn}              { \bfm{n}_m }                      
\newcommand{\ppp}              { \bfm{p} }       
\newcommand{\sig}             { \boldsymbol{\sigma}}           
\newcommand{\vareps}             { \boldsymbol{\varepsilon} }     
\newcommand{\isdef}           { \overset{\text{def}}{=} } 
\newcommand{\ds}              { \displaystyle }   

\newcommand{\uuH}              { \text{H} \bfm{u}  }                  
\newcommand{\ff}             { \bfm{f}  }      %
\newcommand{\GI}          { \Gamma_{\text{I}} }            
\newcommand{\GO}          { \Gamma_{\text{O}} }            
\newcommand{\ccc}             { \bfm{c}  }      

\newcommand{\SHR}            { H^r(\Omega) }                       
\newcommand{\SHS}            { H^s(\Omega) }                       

\theoremstyle{plain}
\theoremheaderfont{\bfseries \upshape}
\newtheorem{thm}{Theorem}[section]
\newtheorem{lem}{Lemma}[section]
\newtheorem{rem}{Remark}[section]

\newtheorem{prp}{Proposition}[section]

\usepackage[margin=1in]{geometry}

\begin{document}

\begin{frontmatter}
 \title{An Adaptive Stable Space-Time FE Method for the Shallow Water Equations}

\author{Eirik Valseth\corref{cor1}}
\ead{Eirik@utexas.edu}

\author{Clint Dawson}
\ead{Clint@oden.utexas.edu}

 \cortext[cor1]{Corresponding author}

 \address{Oden Institute for Computational Engineering and Sciences, The University of 
 Texas at Austin, Austin, TX 78712, USA}


\biboptions{sort&compress}

\begin{abstract}
We consider the finite element (FE) approximation of the shallow water equations (SWE) by considering 
discretizations in which both space and time are established using an unconditionally stable FE method. 
Particularly, we consider the automatic variationally stable FE (AVS-FE) method, a type of 
discontinuous Petrov-Galerkin (DPG) method. 
The philosophy of the DPG method allows us to break the test space and achieve 
unconditionally stable FE approximations as well as accurate \emph{a posteriori} error estimators 
 upon solution of a saddle point system of equations. 
The resulting error indicators allow us to employ mesh adaptive strategies and perform 
space-time mesh refinements, i.e., local time stepping. 

We derive \emph{a priori} error estimates for the AVS-FE method and linearized SWE and perform 
numerical verifications to confirm corresponding asymptotic convergence behavior. 
In an effort to keep the computational cost low, we consider an alternative 
space-time approach in which the space-time domain is partitioned into finite sized space-time slices.
Hence, we can perform adaptivity on each individual slice to preset error tolerances as needed for a
particular 
application. Numerical verifications comparing the two alternatives indicate the space-time slices 
are superior for simulations over long times, whereas the solutions are indistinguishable for short 
times. 
Multiple numerical verifications show the adaptive mesh refinement capabilities of the AVS-FE method,
as well the application of the method to  commonly applied benchmarks for the SWE.

\end{abstract}

\end{frontmatter}

\section{Introduction}
\label{sec:introduction}

The shallow water equations govern the flow of water in domains in which the 
characteristic wavelength horizontally is significantly larger than the depth of water. 
A very important application of the SWE is in the modeling of events such as storm surges resulting
from hurricanes. Thus, the importance of accurate numerical solution techniques should therefore 
be clear as the repercussions of such events can be vast. 
The SWE are models for the Navier-Stokes equations in which the direction of the depth has been 
integrated from the sea floor to the free surface of the water and application of the corresponding boundary conditions.
The resulting equations are continuity and momentum equations for the water
surface elevation and depth-averaged horizontal velocities, respectively. 
Generally, the domains if interest in the application of the SWE are irregular and the resulting 
computational meshes need to be unstructured, thereby making FE methods well suited for the numerical 
approximation of the SWE, see, e.g.,~\cite{peraire1986shallow,kawahara1982selective} for early
examples.

For flow regimes resulting from hurricanes, the motion of the water, i.e., convection is  
the driving mechanism of transport.
The domination of convective transport over diffusive transport leads to stability issues in classical
Galerkin FE methods.
This issue is well known to be resolved when the element size in the
FE mesh is adequately refined near the interior or boundary layers induced by this phenomenon. 
However, the computational cost and the required mesh generation efforts on a case-by-case basis
is generally prohibitive. Furthermore, the smallest element size dictates the corresponding 
time step size further increasing the computational cost.
Lynch and Gray developed the  wave continuity equation as a surrogate for the SWE
in~\cite{lynch1979wave} which result in more stable FE approximations as it changes the enforcement of
 the continuity
equation by introducing second order time derivatives thereby circumventing the hyperbolic nature 
of the SWE. 
The wave continuity equation is employed in the advanced circulation model of
Luettich \emph{et. al}~\cite{luettich1992adcirc}, 
which is widely used and has been developed to encompass a wide range of features important in 
hurricane storm surge modeling.  

Discontinuous Galerkin (DG) methods as well as their hybridized 
versions~\cite{samii2019comparison,gassner2016well,wintermeyer2017entropy} are popular 
FE methods for the SWE.
Important reasons for their popularity include high accuracy, their mass conservation property, ease 
of establishing conditionally
stable FE discretizations, and local $p$-adaptivity~\cite{kubatko2009dynamic}.
Recently, high order entropy stable DG methods for curved FE meshes have been introduced for the SWE 
by Wu \emph{et al.} in~\cite{wu2020high} which satisfies discrete conservation of entropy. The 
entropy stable method takes advantage of summation-by-parts operators to increase 
computational efficiency compared to traditional DG methods demonstrated through numerical
verifications.  
The development of hybridized DG (HDG) methods has also reduced the computational 
cost of solving the global system of equations 
compared to classical DG methods~\cite{samii2019comparison}. 
Coupled Galerkin and DG methods have also successfully been 
developed to maximize efficiency and accuracy by Dawson and Proft 
in~\cite{dawson2002discontinuous,dawson2003discontinuous}.
Least squares FE methods (LSFEMs)~\cite{bochevLeastSquares} have also been successfully applied to the
SWE by Starke in~\cite{starke2005first} and Liang and Hsu in~\cite{liang2009least}. Both these LSFEMs
take advantage of the stability property of LSFEMs spatially and the authors present several numerical
verifications.

All the aforementioned FE methods used in shallow water systems for the SWE or its surrogate 
wave continuity equation employ a method of lines approach in the temporal discretization. 
Thus, spatial and temporal computations are decoupled, where finite elements are employed in space 
and time stepping schemes such as finite difference methods are employed in time.  
It is significantly less common to employ space-time FE methods, i.e., using FE discretizations of 
both space and time. The main reasons are likely the increased computational cost of such methods, as
well as the inherently unstable numerical nature of classical FE methods for first order 
partial derivatives. 
However, some examples of space-time FE methods for the SWE do exist in 
literature~\cite{ribeiro1998space,ribeiro1970finite,takase2010space,arabshahi2016space}, where 
conditional 
discrete stability is ensured in space and time by an upwinding argument. 
In modern multi processor computers and supercomputers, the additional cost of space-time FE methods 
can be justified as the functional framework of FE leads to readily available 
\emph{a priori} error bounds in addition to \emph{a posteriori} error estimation 
techniques. Hence, space-time FE methods can take advantage of adaptive mesh refinement strategies 
to ensure maximum computational efficiency and accuracy.

The AVS-FE method~\cite{CaloRomkesValseth2018} is an unconditionally stable FE method that falls
into the category of discontinuous Petrov-Galerkin (DPG) methods, introduced by Demkowicz
and Gopalakrishnan~\cite{Demkowicz4, Demkowicz2, Demkowicz3, Demkowicz5, Demkowicz6}. 
The unconditional stability of the method is ensured by a particular choice of test functions 
which is defined by a Riesz representation problem that realize the supremum in the
 \emph{inf-sup} condition. 
Hence, the AVS-FE method is a Petrov-Galerkin method in which the weak formulation is such that the
 trial space is a globally continuous Hilbert space, whereas the test space is a broken Hilbert space. 
Hence, the test functions are square integrable functions  that are of higher order 
regularity on each element, e.g., $H^1$ or $H(\text{div})$. The corresponding FE discretization of this 
test space is achieved by by computing on-the-fly optimal test functions 
in the spirit of the DPG method.
This method remains attractive in particular due to its unconditional stability property, regardless 
of the differential operator. 
Additionally, according to the philosophy of the DPG method, the AVS-FE method establishes an 
equivalent saddle point system which yields both the AVS-FE approximation, as well as 
an "error representation function". This function leads to \emph{a posteriori} error estimates of the 
numerical approximation error in terms of the energy norm induced by the sesquilinear form.

In this paper, we build upon the work of~\cite{valseth2020unconditionally}, 
where the AVS-FE method was employed to establish space-time FE approximations for the 
Korteweg-de Vries equation. 
Following this introduction, we introduce the model initial boundary value problem (IBVP) as well as
 notations and conventions in Section~\ref{sec:model_problem}.
 In Section~\ref{sec:avs-discretization}, we present the equivalent AVS-FE weak 
formulation and its analysis for the SWE IBVP as well as the concepts
of Carstensen \emph{et al.}~\cite{carstensen2018nonlinear} to establish nonlinear 
AVS-FE approximations.
\emph{A priori} error estimates are introduced in Section~\ref{sec:avs-estimates}.
In Section~\ref{sec:verifications}, we  perform multiple  numerical verifications for the SWE  
presenting numerical asymptotic convergence properties as well $h-$adaptive mesh refinements.
Finally, we conclude with remarks on the results and future works in Section~\ref{sec:conclusions}.

\section{The Shallow Water Equations and the AVS-FE Method}  
\label{sec:model_and_conv}

We consider a simplified form of the SWE to establish an 
unconditionally stable space-time FE method for the SWE. 
After establishing the SWE model problem, we introduce the corresponding AVS-FE weak formulation and 
discretization.

\subsection{Model Problem: The Shallow Water Equations}
\label{sec:model_problem}

\begin{figure}[h!]
\centering
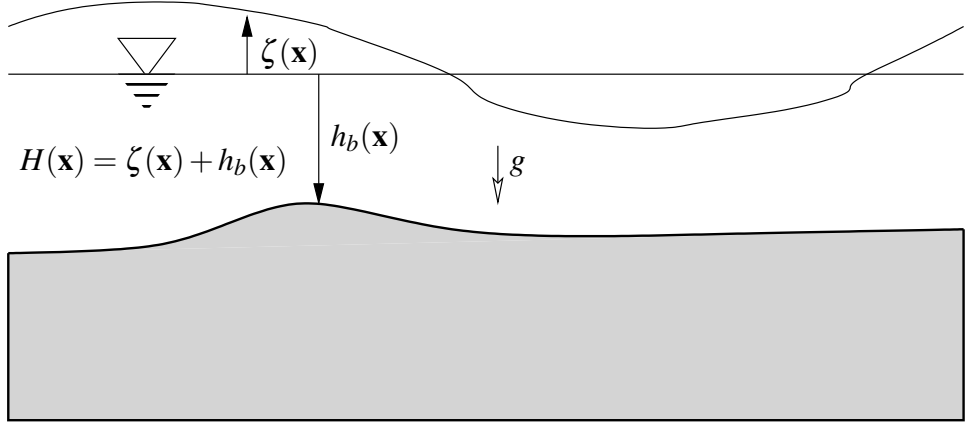 
\caption{Elevation and bathymetry overview.}
\label{fig:domain_bath}
\end{figure}

The derivation of the SWE from the three-dimensional incompressible Navier–Stokes equations
is performed 
under the assumptions of a long horizontal wavelength and a hydrostatic pressure distribution, see, e.g.~\cite{falconer1993introduction}. 
Let $\Omega\subset \mathbb{R}^2$ be a bounded open domain with a Lipschitz boundary  $\partial \Omega$
which is  
partitioned into two segments $\GI$ and $\GO$, such 
that $\partial \Omega = \overline{ \GI \cup \GO}$. Also, let $\vn$ be the outwards unit normal vector 
to the global boundary, and identify the boundary segments as
 $\GI = \{\xx \in \partial \Omega \, : \, \uu \cdot \vn < 0  \}$ and
  $\GO = \{\xx \in \partial \Omega \, : \, \uu \cdot \vn \ge 0  \}$ as in and outflow boundaries, 
respectively.
Finally, define the
temporal domain $t \in (0,T) \subset \mathbb{R}^+_0$.
Hence, we consider the following version of the viscous two dimensional 
SWE~\cite{dawson2003discontinuous,dawson2002discontinuous}:
\begin{subequations}  \label{eq:ShallowWE}
\begin{align}  
\ds \frac{\partial \zeta}{\partial t} + \Nabla \cdot (\uuH)  = 0, \label{eq:continuity}\\ 
\ds \frac{\partial \uu}{\partial t} + \uu \cdot (\Nabla \uu) + \tau_{bf} \, \uu + g \Nabla \zeta 
- \mu \Delta \uu  = \ff, \label{eq:momentum} 
\end{align}
\end{subequations}
where $H = \zeta + h_b$, $h_b = h_b(\xx)$ is the bathymetry of the bottom surface (see Figure~\ref{fig:domain_bath}),
the unknowns $\zeta = \zeta(\xx,t)$ and $\uu=\uu(\xx,t)=\{u_x(\xx,t),u_y(\xx,t) \}^{\text{T}} $ represent depth
averaged elevation and velocity, respectively, $g = 9.81 m/s^2$ the constant of gravitational acceleration, $\mu$ the 
depth averaged turbulent viscosity, $\tau_{bf}$ the bottom friction factor, and $\ff$ represent body
forces. %
The system~\eqref{eq:ShallowWE} consists of two partial differential equations (PDEs), the continuity
equation~\eqref{eq:continuity} for the depth averaged water column elevation, and the momentum
equation~\eqref{eq:momentum} governing the depth averaged velocity.
The bottom friction factor $\tau_b$ is a source of another potential 
nonlinearity as multiple friction models depend on both velocity and water depth.

To establish an IBVP of the SWE with a unique 
solution, a proper combination of boundary conditions (BCs) and initial conditions (ICs) is needed.
Since the SWE are derived from the 
Navier-Stokes equations, these conditions are often not trivial to establish due to the chaotic nature of these 
PDEs. 
In this paper, we seek to establish stable FE approximations of the SWE, therefore, for simplicity, 
 we only consider cases in which the resulting IBVP has a unique solution~\cite{praagman1979numerical}.
To this end, we consider the following boundary and initial conditions:
\begin{equation} \label{eq:SWE_IBCs}
\begin{array}{ll}
\ds \zeta = \hat{\zeta} \text{ on } \GI,&\\  
\ds \uu = \hat{\uu} \text{ on }  \partial \Omega, &\\ 
\ds \zeta = \zeta_0 \text{ on } \Omega, &\\ 
\uu = \uu_0 \text{ on }  \Omega. &\\  
 \end{array}
\end{equation}

Combining the PDE~\eqref{eq:ShallowWE} and the conditions~\eqref{eq:SWE_IBCs} gives the SWE IBVP:
\begin{equation} \label{eq:SWE_IBVP}
\boxed{
\begin{array}{l}
\text{Find }  (\zeta,\uu)  \text{ such that:}    
\\[0.05in] 
\qquad 
\begin{array}{rrl}
\\[-10pt] 
\ds \frac{\partial \zeta}{\partial t} + \Nabla \cdot (\uuH)  = 0, \text{ in } 
\Omega \times (0,T), & \\
\ds \frac{\partial \uu}{\partial t} + \uu \cdot (\Nabla \uu) + \tau_{bf} \, \uu + g \Nabla \zeta 
- \mu \Delta \uu  = \ff, \text{ in } \Omega \times (0,T), & \\ \\[-5pt]
\ds \zeta = \hat{\zeta} \text{ on } \GI,&\\  
\ds \uu = \hat{\uu} \text{ on }  \partial \Omega, &\\ 
\ds \zeta = \zeta_0 \text{ on } \Omega, &\\ 
\uu = \uu_0 \text{ on }  \Omega. &\\  
 \end{array}
 \end{array}
}
\end{equation}
%
In the following, we shall use the following notations:
\begin{itemize}

\item inner products between vector
valued functions are denoted with the single dot symbol "$\cdot$", 
and inner products between tensor valued functions are denoted by the colon or double dot 
symbol "$\colon$". 

\item the operator $\Nabla$ is the spatial gradient operator.

\item the operation $\Nabla \cdot \sig$, $\sig$ being a matrix/tensor valued function, corresponds to a 
row-wise application of the gradient operator. 

\item $h_m$ is the diameter  of element $K_m$.

\item $\nn$ is the outwards unit normal vector to element $K_m$


\end{itemize}

\subsection{Weak Formulation and AVS-FE Discretization}  
\label{sec:avs-discretization}

The derivation of a weak formulation for the AVS-FE method can be found in~\cite{CaloRomkesValseth2018}
 and we omit the step-by-step 
derivation found in these for the SWE and mention key points only. The AVS-FE weak formulations are
established by techniques used in mixed FE methods as well as the broken  weak forms associated to
DPG and DG methods.
The first step is the partition of the computational domain into finite elements, i.e., into  
a FE mesh. Since we are taking a space-time approach, we first define the 
space-time domain $\Omega_{\text{T}}$:
\begin{equation}
\notag
  \Omega_{\text{T}} \isdef \Omega \times (0, \text{T} ),
\end{equation}
and the partition $\Ph$ of $\Omega_{\text{T}}$ into elements $K_m$, is such that:
\begin{equation}
\label{eq:domain_part}
  \Omega_{\text{T}} = \text{int} ( \bigcup_{\Kep} \overline{K_m} ), \quad K_m \cap K_n, \quad m \ne n.
\end{equation}
Our goal is to employ classical continuous FE approximating functions as bases for the trial space,
e.g., Lagrange or Raviart-Thomas functions. Hence, it is required to recast the
IBVP~\eqref{eq:SWE_IBVP} into a first order system by introducing a tensor-valued auxiliary 
variable $\sig$:
\begin{equation}
\label{eq:sigvar}
  \sig = \Nabla \uu   \isdef 
  \begin{bmatrix}
       \ds    \frac{\partial u_x}{\partial x} & \ds \frac{\partial u_x}{\partial y} \\
       \ds    \frac{\partial u_y}{\partial x} & \ds \frac{\partial u_y}{\partial y} \\
\end{bmatrix}. 
\end{equation}
Using this variable, we recast the SWE IBVP as:
\begin{equation} \label{eq:SWE_IBVP_first_order}
\boxed{
\begin{array}{l}
\text{Find }  (\zeta,\uu, \sig)  \text{ such that:}    
\\[0.05in] 
\qquad 
\begin{array}{rrl}
\\[-10pt] 
\ds \frac{\partial \zeta}{\partial t} + \Nabla \cdot (\uuH)  = 0, \text{ in } \Omega_{\text{T}}, & \\
\ds \frac{\partial \uu}{\partial t} + \uu \cdot (\Nabla \uu) + \tau_{bf} \, \uu + g \Nabla \zeta 
- \mu \Nabla \cdot \sig   = \ff, \text{ in } \Omega_{\text{T}}, & \\
\ds \sig - \Nabla \uu = \bfm{0}, \text{ in } \Omega_{\text{T}}, & \\ \\[-5pt]
\ds \zeta = \hat{\zeta} \text{ on } \GI,&\\  
\ds \uu = \hat{\uu} \text{ on }  \partial \Omega, &\\ 
\ds \zeta = \zeta_0 \text{ on } \Omega_0, &\\ 
\uu = \uu_0 \text{ on }  \Omega_0. &\\  
 \end{array}
 \end{array}
}
\end{equation}
Hence, in the weak enforcement of~\eqref{eq:SWE_IBVP_first_order}, the required regularities of the trial
variables are $\zeta \in \SHOOT$, $\uu \in \SHOOT^2$, and $\sig \in \SHdivO$, where $\SHOOT$ 
and $\SHdivO$ are the classical $H^1$ and $H(div)$ Hilbert spaces on $\Omega_{\text{T}}$
and $\Omega$, respectively.

The AVS-FE weak formulation is the established by an element-wise weak enforcement 
of~\eqref{eq:SWE_IBVP_first_order}, application of Green's identity to all spatial derivatives 
with the exception of the convective term $\uu \cdot (\Nabla \uu)$. 
The BCs are all enforced in a weak sense, and we incorporate the initial conditions strongly 
in the definitions of the trial space. Thus, the AVS-FE weak formulation is:
\begin{equation} \label{eq:weak_form}
\boxed{
\begin{array}{ll}
\text{Find } (\zeta, \uu, \sig) \in \UUUT \text{ such that: }
\\[0.05in]
   B((\zeta, \uu, \sig);(v,\www,\ppp) ) = F(v,\www,\ppp), \quad \forall (v,\www,\ppp) \in \VV,
 \end{array}}
\end{equation}
where the sesquilinear  form, $B:\UUUT\times\VV\longrightarrow \mathbb{R}$, and linear 
functional, $F:\VV\longrightarrow \mathbb{R}$ are defined:
\begin{equation} \label{eq:B_and_F}
\begin{array}{l}
\hspace{-0.5cm}B((\zeta, \uu, \sig),(v,\www,\ppp) ) \isdef
\ds \summa{\Kep}{} \int_{K_m} \biggl\{   \frac{\partial \zeta}{\partial t} v_m - \Nabla v_m 
\cdot (\uuH)  \,\ds \, 
\\[0.15in] \ds + [ \frac{\partial \uu}{\partial t} + \uu \cdot (\Nabla \uu) + \tau \, \uu ] 
\cdot \www_m  - g \zeta (\Nabla \cdot \www_m) + \mu   \sig : \Nabla \www_m \,  \\[0.1in] 
\ds      \,+ \,  \sig : \ppp_m + \uu \cdot ( \Nabla \cdot \ppp_m )
 \biggr\}\dx \,    
 \ds + \oint_{\dKm} \biggl\{  v_m( \uuH \cdot \nn )  \ds - \mu (\sig \nn) \cdot \www_m  \biggr\}\dss
  \\ \ds + \oint_{\dKm \setminus \Gamma_{out}} g \zeta (\www_m \cdot \nn) \dss,
 \\[0.15in]
\hspace{-0.5cm}F(v,\www,\ppp) \ds =\summa{\Kep}{} \int_{K_m} \ff \cdot \vv_m \dx  + \ds \oint_{\dKm
 \cap \partial \Omega}  (\hat{\uu} \ppp_m ) \cdot \nn  \dss  - \ds \oint_{\dKm \cap  \Gamma_{in}}  
g \hat{\zeta} (\www_m \cdot \nn) \dss,
 \end{array}
\end{equation}
and the function spaces $\UUUT$ and $\VV$ are defined:
\begin{equation}
\label{eq:function_spaces}
\begin{array}{c}
\UUUT \isdef \biggl\{ (\zeta, \uu, \sig) \in  \SHOOT^3\times \SHdivO:
 \uu_{|\Omega_0 } =\uu_0, \zeta_{|\Omega_0 } =\zeta_0 \ \biggr\},
\\[0.15in]
\VV \isdef \SHOP^3 \times \SHdivPh,
\end{array}
\end{equation}
where broken spaces are defined:
\begin{equation}
\label{eq:broken_h1_space}
\SHOP \isdef \biggl\{ v\in\SLTOT: \quad v_m \in \SHOK, \; \forall \Kep\biggr\},
\end{equation}
\begin{equation}
\label{eq:broken_hdiv_space}
\SHdivPh \isdef \biggl\{ \vv \in\SLTO^4 : \quad \vv_m \in \SHdivK, \; \forall \Kep\biggr\}.
\end{equation}
We also define  norms on these spaces, 
$\norm{\cdot}{\UUU}:  \UUU \!\! \longrightarrow \!\! [0,\infty)$ and
$\norm{\cdot}{\VV}: \VV\! \! \longrightarrow\! \! [0,\infty)$ as:
\begin{equation}
\label{eq:norms}
\begin{array}{l}
\hspace{-0.5cm}\ds \norm{(\zeta,\uu,\sig)}{\UUUT} \isdef \sqrt{\int_{\Omega} \biggl[   \Nabla \zeta
 \cdot \Nabla \zeta + \zeta^2 + \Nabla \uu : \Nabla \uu + \uu \cdot \uu   + (\Nabla \cdot \sig)^2+\sig
  : \sig  \biggr] \dx },
\\[0.2in]
\hspace{-2.5cm} \ds   \norm{(v,\www,\ppp)}{\VV} \isdef \sqrt{\summa{\Kep}{}\int_{K_m} \biggl[  h_m^2
 \Nabla v_m \cdot \Nabla \vv_m + v_m \cdot v_m +h_m^2 \Nabla \www_m : \Nabla \www_m   + \www_m \cdot
 \www_m +  h_m^2 (\Nabla \cdot \ppp_m)^2+\ppp_m : \ppp_m \biggr] \dx }.
 \end{array}
\end{equation}
%
%
%
Note that the norm $\norm{\cdot}{\VV}$ is equivalent to the  classical broken norm:
\begin{equation}
\label{eq:norm2}
\begin{array}{l}
\hspace{-0.5cm} \ds   \norm{(v,\www,\ppp)}{V} \isdef \sqrt{\summa{\Kep}{}\int_{K_m} \biggl[ 
 \Nabla v_m \cdot \Nabla \vv_m + v_m \cdot v_m + \Nabla \www_m : \Nabla \www_m  + \www_m \cdot
 \www_m +  (\Nabla \cdot \ppp_m)^2+\ppp_m : \ppp_m \biggr] \dx }.
 \end{array}
\end{equation}
The weak 
formulation~\eqref{eq:weak_form} is closely related to the weak formulations of the DPG method with
ultraweak forms, differing in our choice of trial space. 

To simplify the analyses of the AVS-FE weak formulation, we introduce an equivalent norm on 
the trial space, 
the energy norm $\norm{\cdot}{\text{B}}: \UUUT\longrightarrow [0,\infty)$:
\begin{equation}
\label{eq:energy_norm}
\norm{(\zeta,\uu,\sig)}{\text{B}} \isdef \supp{(\vv,\www,\ppp)\in \VV\setminus \{0\}} 
     \frac{|B((\zeta, \uu, \sig);(v,\www,\ppp) )|}{\quad\norm{(\vv,\www,\ppp)}{\VV}}.
\end{equation}
Before proceeding to the FE discretization of the weak formulation~\eqref{eq:weak_form}, we point out 
that the weak formulation is well posed as it satisfies
the \emph{inf-sup} and continuity conditions in terms 
of the energy norm with an \emph{inf-sup} constant of unity if we consider a linearized version of the
SWE for simplicity. Thus, the linearized form of this  formulation 
satisfies the required conditions of the Babu{\v{s}}ka Lax-Milgram Theorem~\cite{babuvska197finite},
for an in-depth discussion on broken spaces and variational forms we refer 
to~\cite{carstensen2016breaking}. 
We also define the optimal test functions 
$(\hat{e},\hat{\vareps},\hat{\bfm{E}}) \in \VV$, for each$(\zeta, \uu, \sig) \in \UUUT$ as the
 solution of the Riesz representation problem:
\begin{equation}
\label{eq:riesz_problem}
\begin{array}{rcll}
\ds \left(\, (\hat{e},\hat{\vareps},\hat{\bfm{E}});(\vv,\www) \, \right)_\VV &  \! \! =  \! & 
B((\zeta, \uu, \sig),(v,\www,\ppp) ),& \, \forall (v,\www,\ppp)\in\VV. 
\end{array}
\end{equation}
The Riesz representation problem is well posed with unique solutions due to the inner product in the
left hand side (LHS).
Since~\eqref{eq:riesz_problem} is infinite dimensional, the optimal test functions must be computed 
approximately through an FE approximation. Fortunately, since the test space is broken, the solution 
can be performed element wise, thereby removing the need for a global solve for the optimal test 
functions.
While the definition of the energy norm makes the analysis of AVS-FE weak formulations straightforward, 
it is not computable as it is defined through a supremum.  
Thankfully, the definition of the optimal test functions through the Riesz representation problem
ensures the following norm equivalence:
\begin{equation}
\label{eq:norm_equivalence}
\norm{(\zeta,\uu,\sig)}{\text{B}} = \norm{(\hat{e},\hat{\vareps},\hat{\bfm{E}})}{\VV},
\end{equation}
which is readily available for computations once $(\hat{e},\hat{\vareps},\hat{\bfm{E}})$ is known,
see~\cite{Demkowicz2,Demkowicz5} for details on optimal test functions and proof of the norm equivalence.

To establish FE approximations of~\eqref{eq:weak_form}, we make the standard FE choice of a finite 
dimensional subspace $\UUUhT \subset \UUUT$. Since $\UUUT$ consists of Hilbert spaces, we use 
classical FE basis functions, e.g., $C^0$ polynomials and/or Raviart-Thomas bases.
Here, we make no particular choice other than using conforming approximation spaces.
However, the test space is to be constructed using the DPG philosophy through discrete Riesz
representation problems. 
Thus, the AVS-FE discretization of~\eqref{eq:weak_form} is: 
\begin{equation} \label{eq:discrete_form_SWE}
\boxed{
\begin{array}{ll}
\text{Find } (\zeta^h, \uu^h, \sig^h) \in \UUUhT   \text{ such that:}
\\[0.05in]
 B((\zeta^h, \uu^h, \sig^h);(v^*,\www^*,\ppp^*) ) = F(v^*,\www^*,\ppp^*), \quad \forall 
 (v^*,\www^*,\ppp^*) \in \VVh, 
 \end{array}}
\end{equation}
where 
the finite dimensional test space $\VVh \subset \VV$ is spanned by numerical
approximations of the optimal test functions through the Riesz representation
 problems~\eqref{eq:riesz_problem} computed element-by-element. 
Thus, discrete stability and well-posedness of a linearized version of~\eqref{eq:discrete_form_SWE}
 is ensured due to the 
DPG philosophy (see Remark~\ref{rem:discrete_Stab}).
 Since the underlying differential 
operator is nonlinear, the solution procedure for~\eqref{eq:discrete_form_SWE} must perform 
nonlinear iterations, see the work of Roberts \emph{et al.}~\cite{roberts2015discontinuous} for 
a discussion on DPG method and the Navier-Stokes equations.

Computing optimal test functions on-the-fly to assemble the FE system of linear algebraic equations 
is one interpretation of DPG methods. Another equivalent interpretation is 
that of a mixed or saddle point problem in 
which the Riesz representation problem~\eqref{eq:riesz_problem} is used to define a constraint 
to the weak form~\eqref{eq:weak_form}:
\begin{equation} \label{eq:disc_mixed_SWE}
\boxed{
\begin{array}{ll}
\text{Find } (\zeta^h, \uu^h, \sig^h) \in  \UUUhT, 
(\tilde{e}^h,\tilde{\vareps}^h,\tilde{\bfm{E}}^h) \in V^h(\Ph)  & \hspace{-0.15in} \text{ such that:}
\\[0.05in]
   \quad \left(\, (\tilde{e}^h,\tilde{\vareps}^h,\tilde{\bfm{E}}^h),(v^h,\www^h,\ppp^h) \right)_{\VV} \,
 - B((\zeta^h, \uu^h, \sig^h);(v^h,\www^h,\ppp^h) )& = - F(v^h,\www^h,\ppp^h) \quad  \\  &
\hspace{-.75in}\forall (v^h,\www^h,\ppp^h) \in V^h(\Ph),  \\[0.05in]
  \quad B'((a^h,\bb^h,\ccc^h),(\tilde{e}^h,\tilde{\vareps}^h,\tilde{\bfm{E}}^h) )& =  0, 
  \quad \\ & \hspace{-.75in} \forall \, (a^h,\bb^h,\ccc^h)) \in \UUUhT.
 \end{array}}
\end{equation}
The derivation of this system is based on the weak form~\eqref{eq:weak_form}, the
Riesz representation problem~\eqref{eq:riesz_problem}, and the fact that the AVS-FE is a minimum
residual method, see~\cite{demkowicz2014overview} and~\cite{carstensen2018nonlinear} for details on the
derivation for linear and nonlinear problems, respectively.
$B'$ represents the Gateaux derivative of the sesquilinear form with respect to $(\zeta, \uu, \sig)$
acting on the approximation of the  "error representation 
function" $(\tilde{e},\tilde{\vareps},\tilde{\bfm{E}})$, i.e., it represents a linearization of the 
sesquilinear form.
The mixed problem~\eqref{eq:disc_mixed_SWE} is equivalent to~\eqref{eq:discrete_form_SWE} 
(see Theorem 2.2 of Carstensen \emph{et al.}~\cite{carstensen2018nonlinear}). The linearization
provided by the constraint equation on the Gateaux derivative $B'$ allows us to employ existing 
nonlinear FE solvers to perform iterations that (hopefully) converge to a stationary point.

The error representation function is a Riesz representer of the approximation error 
$(\zeta, \uu, \sig)-(\zeta^h, \uu^h, \sig^h)$ through an analogue of~\eqref{eq:riesz_problem} 
where the right hand 
side is the residual functional.
Thus, due to the norm equivalence~\eqref{eq:norm_equivalence} the norm of the approximate error 
representation function is an \emph{a posteriori} error estimate of the approximation error 
of the energy norm:
\begin{equation}
\label{eq:energynorm_eqv}
\norm{(\zeta-\zeta^h, \uu-\uu^h, \sig - \sig^h)}{\text{B}} \approx 
\norm{(\tilde{e}^h,\tilde{\vareps}^h,\tilde{\bfm{E}})^h}{\VV}.
\end{equation}
This error estimate  has been analyzed in~\cite{Demkowicz2} and its 
and its local restriction to an element is an error indicator:
\begin{equation}
\label{eq:enrr_ind}
\eta = \norm{(\tilde{e}^h,\tilde{\vareps}^h,\tilde{\bfm{E}})^h}{\VVK}.
\end{equation}
This error indicator has been successfully applied to a wide range of problems in both the DPG and 
AVS-FE methods
~\cite{carstensen2016breaking,roberts2015discontinuous,valseth2020stable,valseth2020unconditionally,carstensen2018nonlinear}.
\begin{rem}
The size of the mixed discrete system of linear algebraic equations~\eqref{eq:disc_mixed_SWE}
is larger than the linear system corresponding to~\eqref{eq:discrete_form_SWE} since we do not compute optimal test functions on the fly element-by-element.
An advantage of this mixed form is that its solution immediately provides an a posteriori error
estimate as well as error indicators to be used in mesh adaptive strategies.
\end{rem}
An application of the definition of the Gateaux derivative yields the operator 
$B':\UUUT\times\VV\longrightarrow \mathbb{R}$:
\begin{equation} \label{eq:B_prime_KDV}
\begin{array}{l}
B'((a,\bb,\ccc),(v,\www,\ppp ) \isdef 
\ds \summa{\Kep}{} \int_{K_m} \biggl\{   \frac{\partial a}{\partial t} v_m - \Nabla v_m \cdot ( a 
\uu + \bb H )  \,\ds \, 
+ [ \frac{\partial \bb}{\partial t}  + \tau \, \bb ] \cdot \www_m  \\[0.15in] \ds 
+ \ds  \{ u_x \frac{\partial b_x}{\partial x} + b_x \frac{\partial u_x}{\partial x} + u_y
 \frac{\partial b_x}{\partial y} + b_y \frac{\partial u_x}{\partial y}  ,  u_x \frac{\partial b_y}
 {\partial x} + b_x \frac{\partial u_y}{\partial x} + u_y \frac{\partial b_y}{\partial y} + b_y
 \frac{\partial u_y}{\partial y}   \}^T   \cdot \www_m  \\[0.15in]
 - g a (\Nabla \cdot \www_m) + \mu   \ccc : \Nabla \www_m \,  
\ds      \,+ \,  \ccc : \ppp_m + \bb \cdot ( \Nabla \cdot \ppp_m ) 
 \biggr\}\dx \,    \\[0.1in] 
 \ds + \oint_{\dKm} \biggl\{ ( a \uu + \bb H) \cdot \nn  v_m  \ds - \mu (\ccc \nn) 
 \cdot \www_m  \biggr\}\dss
\ds + \oint_{\dKm \setminus \Gamma_{out}} g a (\www_m \cdot \nn) \dss.
 \end{array}
\end{equation}
\begin{rem} \label{rem:discrete_Stab}
The discrete stability due to optimal test functions is unconditional in the "ideal" case in which 
we can exactly compute these functions. However, since this is not achievable in practical 
computations, we are forced to consider a practical implementation and consider 
an approximation of these functions~\cite{Gopalakrishnan1}. 
Hence, the approximation of the optimal test functions leads to a potential loss of discrete 
stability of the approximation is not sufficiently accurate. In the DPG method, sufficient 
accuracy of the optimal test functions is guaranteed by the existence of (local) Fortin
operators~\cite{boffi2013mixed}. The construction of such operators is studied in great detail in 
~\cite{nagaraj2017construction}, and its analysis was recently further refined 
in~\cite{demkowicz2020construction}. Generally, in DPG methods for second order PDEs, a Fortin 
operator's existence and thus discrete stability is ensured if the local Riesz representation 
problems are solved using polynomials of order $r=p+\Delta p$, where $p$ is the degree of the trial 
space discretization and $\Delta p = d$ the space dimension. However, while this enrichment degree
ensures the existence of the required Fortin operator, numerical evidence suggest that in most cases
$\Delta p = 1$ is typically sufficient~\cite{demkowicz2020construction}. Alternative test 
spaces for the DPG method for singular perturbation problems are investigated 
in~\cite{salazar2019alternative}, even for the case of $\Delta p = 0$.
Another consequence of the existence of Fortin operators is the robustness of the \emph{a posteriori}
error estimate through the error representation function~\eqref{eq:energynorm_eqv}
as shown in~\cite{Demkowicz2}. 

In  the AVS-FE method, numerical evidence suggests that $r=p$ is 
sufficient~\cite{CaloRomkesValseth2018,valseth2020goal} for convection-diffusion PDEs as well as extensive numerical experimentation for the SWE. Since the test functions are sought in a
discontinuous polynomial space, using $r=p$ still results in a larger space than the trial as the 
discontinuous spaces contain additional degrees of freedom. Furthermore, in the limit $h \rightarrow 0$
the space $\VV$ is essentially $L^2$, i.e, any polynomial degree above constants is inherently 
an enrichment of the test space. 
%
\end{rem}

\subsection{Time Slice Approach}
\label{sec:time_slices}

The numerical solution of the space-time discretizations~\eqref{eq:discrete_form_SWE} 
or~\eqref{eq:disc_mixed_SWE} is akin to a FE discretization of a three dimensional problem.
In an effort to keep the global number of degrees of freedom low, we also introduce a "time-slice" 
approach as introduced for transient convection-diffusion in~\cite{ellis2016space,ellis2014space} for
the DPG method. 
Similar approaches were also employed 
in~\cite{ribeiro1998space,ribeiro1970finite,takase2010space,arabshahi2016space} for the SWE for 
stabilized DG and streamlined-upwind Petrov-Galerkin (SUPG) methods.

\begin{figure}[h!]
\centering
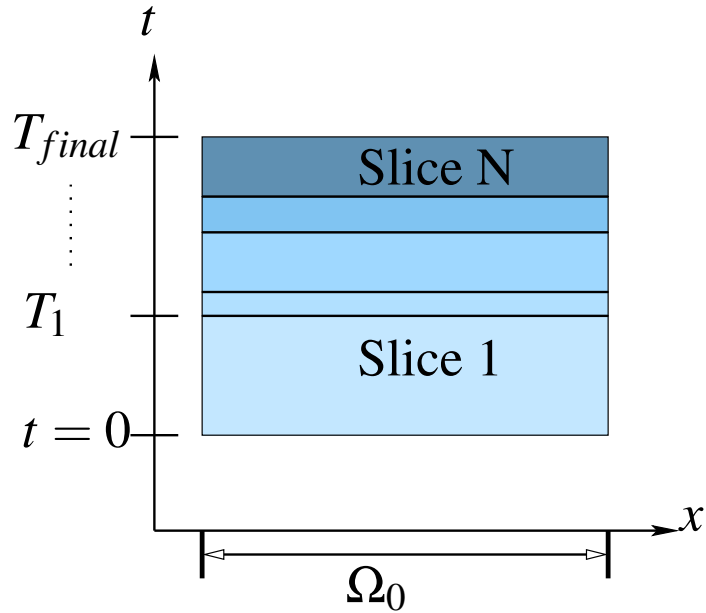 
\caption{Partition of a space-time domain into slices.}
\label{fig:space_time_slices}
\end{figure}

Since solution information is unidirectional in time, we are enabled to consider a partition of the 
temporal part of out domain $\Omega_T$ into "slices", see Figure~\ref{fig:space_time_slices} for an 
example in a spatially one dimensional domain. 
In each slice, starting with the one  intersecting  the initial time boundary $\Omega_0$, we solve the
AVS-FE discrete problem~\eqref{eq:discrete_form_SWE} as described in 
Section~\ref{sec:avs-discretization}.
Each successive slice is then solved by taking the previous solutions $\zeta^h_{prev}, \uu^h_{prev}$
as initial conditions. 
The slice configuration can be established in several ways, several of which are considered for the 
space-time DPG method in~\cite{ellis2016space,ellis2014space}.
Furthermore, the built-in error estimator and error indicators of the AVS-FE method lets us perform 
space-time mesh adaptive refinements on each individual slice. 
This is of particular interest in applications where the solution response varies greatly in time,
thereby allowing us to adaptively refine the mesh to the required resolution in each slice.
This is shown in Figure~\ref{fig:space_time_slices} with variable thickness of the slices in 
the $t$ direction.
In Section~\ref{sec:verifications} we present numerical verifications for the time-slice approach.

\section{A Priori Error Analysis}  
\label{sec:avs-estimates}

Here, we present \emph{a priori} error estimates for the full space-time AVS-FE approximations for
 the SWE in terms 
of appropriate norms of the approximation error. Due to the energy norm~\eqref{eq:energy_norm} and the 
best approximation error  of the AVS-FE method in terms of this norm, the following proofs  
rely on classical bounds in Hilbert spaces for continuous and discontinuous FE approximation functions. 
The presentation is self contained, and researchers considering a method of lines approach with 
AVS-FE method in space while using time stepping schemes can apply these results in their analyses.
In the analysis we assume that the elevation and velocity trial variables are discretized using 
classical $C^0$ continuous
 FE polynomial basis functions, whereas $\sig^h$ is discretized using Raviart-Thomas bases. 
Other conforming alternatives for $\sig^h$ such as Brezzi-Douglas-Marini functions can also be 
employed with minor modifications to the following analyses and we refer to the text by Brezzi 
and Fortin~\cite{BrezziMixed} for details.
We also assume that we are considering a linearized version of the SWE.

First, we introduce the best and quasi-best approximation properties of the AVS-FE method to
facilitate proofs of the \emph{a priori} bounds which are to follow. 
\begin{prp} \label{prp:best_approx} 
Let $(\zeta, \uu, \sig) \in \UUUT$ be the exact solution of the AVS-FE weak 
formulation~\eqref{eq:weak_form} and  $(\zeta^h,\uu^h, \sig^h) \in \UUUhT$ its 
corresponding AVS-FE
approximation from~\eqref{eq:discrete_form_SWE} or~\eqref{eq:disc_mixed_SWE}. Then:
\begin{equation} \label{eq:best_approx}
\norm{(\zeta-\zeta^h,\uu-\uu^h,\sig-\sig^h)}{\rm{B}} \leq  
 \norm{(\zeta-v^h,\uu-\www^h,\sig-\ppp^h)}{\rm{B}},
\end{equation} 
where $(v^h,\www^h,\ppp^h) \in \UUUhT$.
\end{prp} 
We refer to~\cite{Demkowicz5} for a proof of this proposition and note that the proof
relies 
on the definition of the energy norm~\eqref{eq:energy_norm}, the Riesz 
problem~\eqref{eq:riesz_problem}, and Galerkin's orthogonality condition.

 \noindent 
Furthermore, since the energy norm is an equivalent norm on $\UUUT$, we also have the quasi-best 
approximation property:
\begin{equation} \label{eq:Qbest_approx}
\norm{(\zeta-\zeta^h,\uu-\uu^h,\sig-\sig^h)}{\UUUT} \leq   C
 \norm{(\zeta-v^h,\uu-\www^h,\sig-\ppp^h)}{\UUUT},
\end{equation} 
where the constant $C$ is independent of the mesh and depends on norm equivalence constants between 
the energy norm and $\norm{\cdot}{\UUUT}$ as well as the continuity constant of a Fortin 
operator (see Remark~\ref{rem:discrete_Stab}).

Another key component in the following analysis is the convergence of polynomial interpolating
functions. Hence,  there exist  a polynomial interpolation operator 
$\Pi_{hp}$~\cite{babuvska1987hp}:
\begin{equation} \label{eq:Pi_inter}
\ds  \Pi_{hp} \, : \, U \rightarrow  U^{hp}. 
\end{equation} 
Thus, $\Pi_{h}(u)$ represents an interpolant of $u$ consisting of continuous 
polynomials,  then~\cite{oden2012introduction}:
\begin{thm} \label{thm:classical_a_priori}
Let $u \in \SHR$ and $ \Pi_{hp}(u)  \in U^{hp}$ be the interpolant of $u$ \eqref{eq:Pi_inter}. 
Then, there exists $C > 0$
such that the interpolation error can be bounded as follows:
%
%
\begin{equation}
\ds  \norm{u-\Pi_{hp}(u)}{\SHS} \leq  C\, \frac{\ds h^{\,\mu-s}}{\ds p^{\,r-s}} \norm{u}{\SHR},
\end{equation} 
where $h$ is the maximum element diameter, $p$ the minimum polynomial degree of
interpolants in the mesh, $s \leq r$, and $\mu = \rm{ min }$ $(p+1,r)$.
\end{thm}
Second, we have the interpolation operator for Raviart-Thomas~\cite{BrezziMixed} spaces, $\rho_{hp}$:
\begin{equation} \label{eq:rho_inter}
\rho_{hp}  \, : \, Q \rightarrow  Q^{hp}, 
\end{equation} 
Thus, $\rho_{hp}(\qq)$ represents an interpolant of $\qq$ consisting of  polynomials which
\textbf{normal} components are continuous, then \cite{BrezziMixed}:
\begin{thm} \label{thm:classical_a_priori_hdiv}
Let $\qq \in H(\text{\bf{div}},\Omega)$ and $ \rho_{hp}(\qq)  \in Q^{hp}$ be the interpolant of 
$\qq$~\eqref{eq:rho_inter}. Then, there exists $C > 0$
such that the interpolation error can be bounded as follows:
\begin{equation}
\ds \exists \, C > 0 \, : \norm{\qq-\rho_{hp}(\qq)}{\SHdivO} \leq  C\, 
\ds h^{\,n}\ds |\qq|_{H^{n+1} ( \Omega)},
\end{equation} 
where $h$ is the maximum element diameter and $n$ the minimum order of Raviart-Thomas interpolants
 in the mesh.
\end{thm}
The final point we highlight before proceeding with the main results of this section are   
on the convergence properties of interpolants of piecewise discontinuous 
polynomials. In~\cite{riviere1999improved} Rivi\'{e}re \emph{et al.} present 
a result analogous to Theorem~\ref{thm:classical_a_priori} for approximations in broken Hilbert spaces.

The first \emph{a priori} bound we present is in terms of the energy norm. While not exactly
computable, it is natural to present as the energy norm is central to the stability of our method.
Furthermore, we can approximate the error in the energy norm through~\eqref{eq:energynorm_eqv}.
To establish error estimates in terms of the energy norm, we first establish a bound on the Riesz 
representers of the trial functions, i.e., the optimal test functions. 
\begin{lem} \label{lem:energy_bound_riesz} 
Let $(\zeta, \uu, \sig) \in \UUUT$ be the exact solution of the AVS-FE weak 
formulation~\eqref{eq:weak_form} and  $(\zeta^h, \uu^h, \sig^h) \in \UUUhT$ its 
corresponding AVS-FE approximation from~\eqref{eq:discrete_form_SWE} or~\eqref{eq:disc_mixed_SWE}.
Then:
\begin{equation} \label{eq:energy_bound}
\norm{(\zeta-\zeta^h,\uu-\uu^h,\sig-\sig^h)}{\rm{B}} \leq C\, \ds 
 \frac{ h^{\,\mu-1}}{ p_{\uu}^{\,r-1}},
\end{equation} 
where $h$ is the maximum element diameter, $\mu =$ $\rm{min}$ $(p_{\uu}+1,r)$,  $p_{\uu}$ the minimum
 polynomial degree of approximation of $\uu^h$ and $\zeta^h$
in the mesh, and  $r$ the minimum regularity  of the solution components 
$(\hat{e},\hat{\vareps},\hat{\bfm{E}})$ of the distributional PDE underlying the
 Riesz representation problem~\eqref{eq:riesz_problem}.
\end{lem} 
\emph{Proof:} The RHS of~\eqref{eq:best_approx} can be bounded by 
the RHS of~\eqref{eq:best_approx} can be bounded by 
the error in the Riesz representers of the exact and approximate AVS-FE trial functions 
by the energy norm equivalence in~\eqref{eq:norm_equivalence}, and the map induced by the Riesz
representation problem~\eqref{eq:riesz_problem} to yield:
\begin{equation*}  
\norm{(\zeta-\zeta^h,\uu-\uu^h,\sig-\sig^h)}{\rm{B}} \leq  C 
\norm{(\hat{e}-\hat{e}^h,\hat{\vareps}-\hat{\vareps}^h,\hat{\bfm{E}}-\hat{\bfm{E}}^h)}{\VV},
\end{equation*} 
where $(\hat{e},\hat{\vareps},\hat{\bfm{E}}) \in \VV$  are the exact Riesz representers of 
$(\zeta,\uu, \sig)$ through~\eqref{eq:riesz_problem}, and 
$(\hat{e},\hat{\vareps},\hat{\bfm{E}}) \in \VVh$  are the approximate Riesz representers 
of $(\zeta^h, \uu^h, \sig^h)$ through a FE discretization of~\eqref{eq:riesz_problem}. 
Hence, the constant $C$ depends on a Fortin type operator~\cite{nagaraj2017construction} which
measures the degree of loss of stability between the continuous and discrete problems.
The definition of the norm $\norm{\cdot}{\VV}$ in~\eqref{eq:norms} and its equivalence to the norm in~\eqref{eq:norm2} then gives:
\begin{equation*}
\begin{array}{lll}
\norm{(\zeta-\zeta^h,\uu-\uu^h,\sig-\sig^h)}{\rm{B}} \leq C \{ \norm{\hat{e}-\hat{e}^h}{\SHOP} +
\norm{\hat{\vareps}-\hat{\vareps}^h}{\SHOP} + \norm{\hat{\bfm{E}}-\hat{\bfm{E}}^h}{\SHdivP}  \}. 
\end{array}
\end{equation*}
Since $\norm{\cdot}{\SHdivP} \leq \norm{\cdot}{\SHOP}$, we get:
\begin{equation*}
\begin{array}{lll}
\norm{(\zeta-\zeta^h,\uu-\uu^h,\sig-\sig^h)}{\rm{B}} \leq C \{ \norm{\hat{e}-\hat{e}^h}{\SHOP} +
\norm{\hat{\vareps}-\hat{\vareps}^h}{\SHOP} + \norm{\hat{\bfm{E}}-\hat{\bfm{E}}^h}{\SHOP}  \}.
 \end{array}
\end{equation*}
Now, we pick $(\hat{e}^h,\hat{\vareps}^h,\hat{\bfm{E}}^h)$ to be polynomial interpolants and apply the
bounds in~\cite{riviere1999improved} for the discontinuous polynomials to get:
\begin{equation*}
\begin{array}{lll}
\norm{(\zeta-\zeta^h,\uu-\uu^h,\sig-\sig^h)}{\rm{B}} \leq C \{ \ds 
 \frac{ h^{\,\mu_{\hat{e}}-1}}{ p_{\hat{e}}^{\,r_{\hat{e}}-1}} + \ds  
\frac{ h^{\,\mu_{\hat{\vareps}}-1}}{ p_{\hat{\vareps}}^{\,r_{\hat{e}}-1}} +
 \ds  \frac{ h^{\,\mu_{\hat{e}}-1}}{ p_{\hat{\bfm{E}}}^{\,r_{\hat{\bfm{E}}}-1}}  \},
 \end{array}
\end{equation*}
where $\mu_i =$ $\rm{min}$ $(p_{i}+1,r_i)$, $r_i$ the regularity of the solution of the underlying 
distributional PDE and $i$ denote 
the components $(\hat{e},\hat{\vareps},\hat{\bfm{E}})$. Finally, we complete the proof by noting that
we that we pick the same degree 
$p$ for all variables and that the term with the smallest $r$ will dominate the error. Hence, we 
get the desired bound~\eqref{eq:energy_bound}.
\noindent ~\qed
\newline \noindent

Second, we establish a bounds in terms of the Sobolev Norm $\norm{\cdot}{\UUUT}$.
\begin{lem} \label{lem:sobolev_error} 
Let $(\zeta, \uu, \sig) \in \UUUT$ be the exact solution of the AVS-FE weak formulation
\eqref{eq:weak_form} and  $(\zeta^h, \uu^h, \sig^h) \in \UUUhT$ its corresponding AVS-FE 
approximation from~\eqref{eq:discrete_form_SWE} or~\eqref{eq:disc_mixed_SWE}. Then:
\begin{equation} \label{eq:sobolev_bound}
\ds \norm{(\zeta-\zeta^h,\uu-\uu^h,\sig-\sig^h)}{\UUUT} \leq C\, \{ \ds  
\frac{ h^{\,\mu-1}}{ p_{\uu}^{\,r-1}} + h^n \},
\end{equation} 
where $h$ is the maximum element diameter, $\mu =$ $\rm{min}$ $(p+1,r)$,  $p$ the minimum polynomial
 degree of approximation of $\uu^h$ and $\zeta^h$
in the mesh,  $r$ the minimum regularity  of the solution components $\uu$ and $\zeta$ of the
underlying distributional SWE PDEs, and $n = p_{\uu}-1$ the minimum order of the Raviart-Thomas elements
in the mesh.
\end{lem} 
\emph{Proof:} By the quasi-best approximation property~\eqref{eq:Qbest_approx}:
\begin{equation*}
\begin{array}{lll}
\ds \norm{(\zeta-\zeta^h,\uu-\uu^h,\sig-\sig^h)}{\UUUT}  \leq C \,
 \norm{(\zeta-v^h,\uu-\www^h,\sig-\ppp^h)}{\UUUT},
 \end{array}
\end{equation*}
the definition of the norm on $\UUUT$~\eqref{eq:norms} leads to :
\begin{equation*}
\begin{array}{lll}
\ds \norm{(\zeta-\zeta^h,\uu-\uu^h,\sig-\sig^h)}{\UUUT} & \leq C \, \{ \ds  \norm{\zeta-v^h}{\SHOOT} + 
\norm{\uu-\www^h}{\SHOOT} + \norm{\sig-\ppp^h}{\SHdivO} \}, 
 \end{array}
\end{equation*}
by choosing functions for $v^h,\www^h$ that are polynomial interpolants as 
in~\eqref{eq:Pi_inter}, and Raviart-Thomas functions for $\ppp^h$~\eqref{eq:rho_inter}.
An application of Theorem~\ref{thm:classical_a_priori} with $s = 1$ and 
Theorem~\ref{thm:classical_a_priori_hdiv} then gives:
\begin{equation*}
\begin{array}{lll}
\ds \norm{(\zeta-\zeta^h,\uu-\uu^h,\sig-\sig^h)}{\UUUT} & \leq C \, \{\ds 
 \frac{ h^{\,\mu_{\zeta}-1}}{ p_{\zeta}^{\,r-1}} 
 +\frac{ h^{\,\mu_{\uu}-1}}{ p_{\uu}^{\,r-1}} + h^n  \}.
 \end{array}
\end{equation*}
In the AVS-FE method, we always pick $p_{\uu} = p_{\zeta}$, then
considering only the largest of the fractions, the bound~\eqref{eq:sobolev_bound} follows
and the proof is complete.
\noindent ~\qed
\newline \noindent


Last, we have error bounds in weaker norms, in this case $L^2$ norms, 
we apply the Aubin-Nitsche lift~\cite{aubin1987analyse,nitsche1972dirichlet} 
for the errors in $\zeta-\zeta^h$ and $\uu-\uu^h$. Whereas for the stress variable $\sig-\sig^h$, 
we rely on established bounds for Raviart-Thomas approximations found in~\citep{BrezziMixed}.
Thus, the following can be established:
\begin{prp} \label{prp:rate_l2}
Let $(\zeta, \uu, \sig) \in \UUUT$ be the exact solution of the AVS-FE weak formulation
\eqref{eq:weak_form} and  $(\zeta^h, \uu^h, \sig^h) \in \UUUhT$ its corresponding AVS-FE 
approximation from~\eqref{eq:discrete_form_SWE} or~\eqref{eq:disc_mixed_SWE}. Then:
\begin{equation} \label{eq:U_rates_l2}
\ds \norm{(\zeta-\zeta^h,\uu-\uu^h,\sig-\sig^h)}{\SLTOT} \leq C\, \{ \ds  
\frac{ h^{\,\mu_{\zeta} -1}}{ p^{\,r_{\zeta}-1}} + \frac{ h^{\,\mu_{\uu}  }}{ p^{\,r_{\uu}}} +  h^{n+1} \},
\end{equation} 
where $h$ is the maximum element diameter, $\mu_{\zeta} =$ $\rm{min}$ $(p_{\zeta}+1,r_{\zeta})$,  $p_{\zeta}$ the minimum polynomial
 degree of approximation of $\zeta^h$
in the mesh,  $r_{\zeta}$ the minimum regularity of the $\zeta$ of the
underlying distributional SWE PDEs, , $\mu_{\uu} =$ $\rm{min}$ $(p_{\uu}+1,r_{\uu})$,  $p_{\uu}$ the minimum polynomial
 degree of approximation of $\uu^h$
in the mesh,  $r_{\uu}$ the minimum regularity of the $\uu$ of the
underlying distributional SWE PDEs, 
 and $n = p_{\uu}-1$ the minimum order of the Raviart-Thomas elements
in the mesh.
where $h$ is the maximum element diameter.
\end{prp} 
\noindent $\square$  

 \noindent 
The individual terms corresponding to the error in the velocity and $\sig$ are optimal in sense 
that it will result in convergence of the same order as their underlying polynomial interpolants.
However, inspection of the term corresponding to the elevation reveals that it is sub optimal 
by one order. The reason is the regularity of a dual solution corresponding to this 
variable with an $L^2$ source is an order lower than the regularity of the dual solution 
corresponding to $\uu$ since the SWE contain derivatives of $\zeta$ of first order only. 
\begin{rem}
To conclude this section, we make an important remark on the AVS-FE approximations for the 
variable $\sig^h$. For convex domains $\Omega_T$ and smooth sources $\ff$, the regularity of $\sig$ is 
higher than the $\SHdivO$ dictated by the weak form~\eqref{eq:weak_form}, i.e., $\SHOO$. Thus, 
it is appropriate to use $C^0$ continuous basis functions as advocated 
in~\cite{CaloRomkesValseth2018}. As the application of the SWE often in complex 
coastal domains which are highly non regular, users of this method should use $\SHdivO$ conforming 
approximations such as Raviart-Thomas elements as a general rule-of-thumb in coastal domains. 
As in~\cite{BrezziMixed}, the minimum order of the Raviart-Thomas elements is one order below the 
polynomial order for the approximate velocity. Hence, increasing the order of these elements above this 
does not increase the accuracy as the bounds on $\sig-\sig^h$ depend on the order of the velocity 
approximations.    
\end{rem}

\section{Numerical Verifications}  
\label{sec:verifications}

In this section, we present numerical verifications for the AVS-FE method and the SWE. First, we
consider several academic problems to ascertain convergence behavior under both uniform and 
adaptive mesh refinements. Finally, we consider a series of common physical benchmark
 problems for shallow water models from literature. 
While the \emph{a priori} error estimates are established based on a linearized SWE, we 
consider only the full nonlinear case here as a linear SWE is generally not appropriate in physical
 applications of interest.
For all verifications, we solve the saddle point system~\eqref{eq:disc_mixed_SWE}.
All numerical experiments presented here are performed using the FE solver
FEniCS~\cite{alnaes2015fenics}.
To converge to the nonlinear SWE solutions, we employ the Portable, Extensible Toolkit for Scientific Computation 
(PETSc) library Scalable Nonlinear Equations Solvers (SNES)~\cite{abhyankar2018petsc,petsc-user-ref} within FEniCS. To minimize the errors of these iterations, we set the nonlinear convergence tolerance 
to the value $10^{-14}$. Based on our experience this typically leads to a converged solution in less than 6 nonlinear iterations.

\subsection{Numerical Asymptotic Convergence Studies}
\label{sec:convergence}

To establish the convergence properties of the AVS-FE method applied to the SWE, we consider a case
where the bathymetry is assumed to be constant equal to zero, i.e., $H=\zeta$.
We consider a case in which the exact elevation is:
\begin{equation} \label{eq:zeta_exact_ex1}
\begin{array}{l}
\ds \zeta_{ex} = \text{cos} (x + y - t),
 \end{array}
\end{equation}
whereas both components of the velocity vector $\uu_{ex}$ are:
\begin{equation} \label{eq:u_exact_ex1}
\begin{array}{l}
\ds u_{ex} = \text{sin} (x + y + t).
 \end{array}
\end{equation}
The other parameters we pick are $\mu = 10^{-5}$ and $\tau_{bf} = 1$, the spatial domain is the
 unit square 
and the temporal domain is from $t=0s$ to $t=0.5s$. 
The exact solutions are then used to establish proper boundary conditions and source terms for the SWE.
In the FE approximation, we employ classical 
Lagrange bases for $\zeta^h$, $\uu^h$, and $\sig^h$ since in this case the smooth solutions make 
$C^0$ regularity a valid choice for discretization of $\SHdivO$. Since we 
solve~\eqref{eq:disc_mixed_SWE}, we also need to specify functions spaces for the error representation
function, in this case we pick discontinuous Lagrange polynomials for all components of this space of
 the same polynomial degree as the velocity.
The initial mesh consists of six
tetrahedron elements to which we perform uniform mesh refinements to ascertain the asymptotic 
convergence rates of applicable norms of the numerical approximation errors.
\begin{figure}[h!]
\subfigure[ \label{fig:SWE_uniform_p2_l2} $\SLTOT$ and Energy error norms. ]{\centering
 \includegraphics[width=0.5\textwidth]{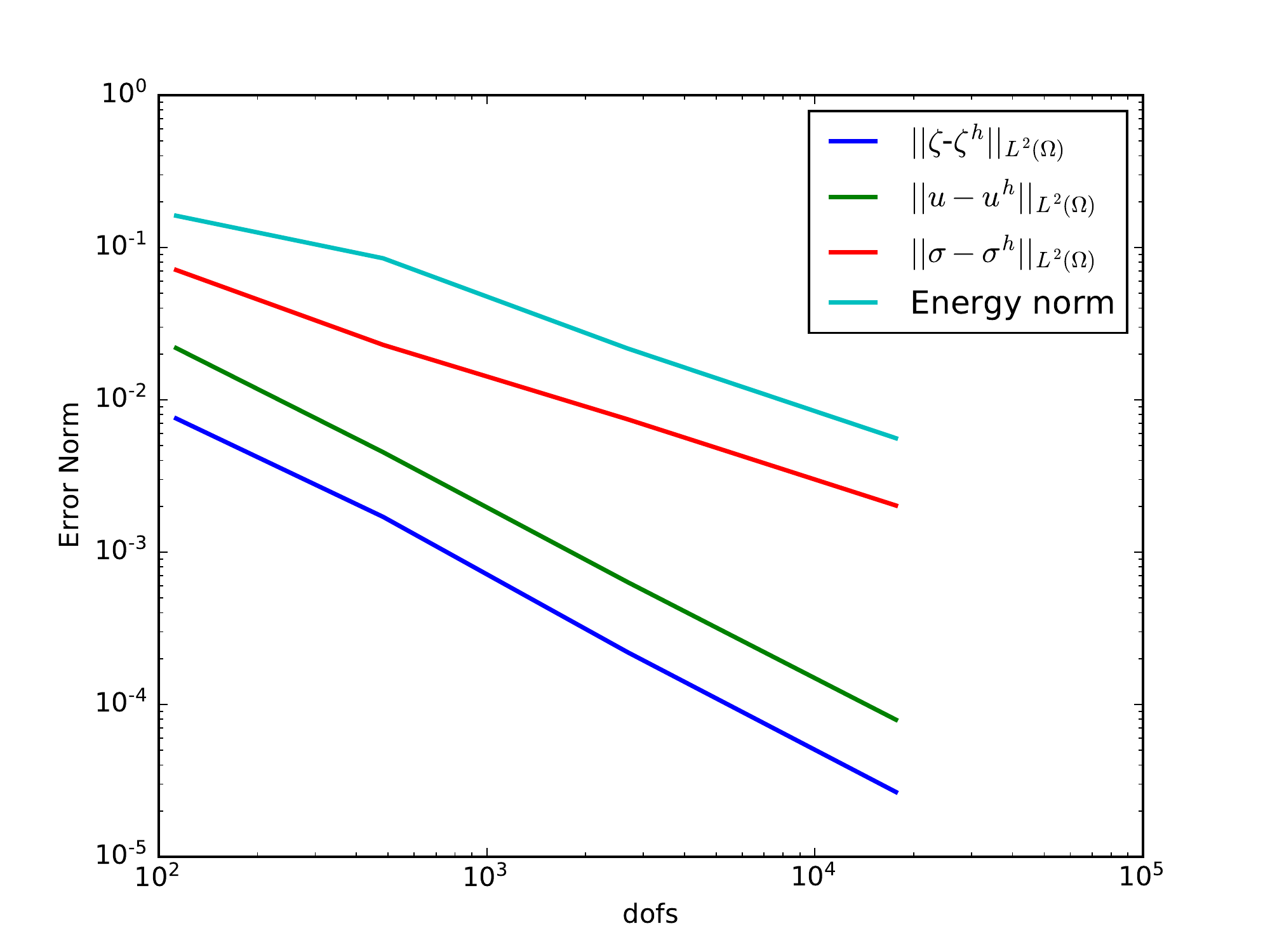}}
  \subfigure[ \label{fig:SWE_uniform_p2_h1} $\SHOOT$, $\SHdivO$ and Energy error norms.]{\centering
 \includegraphics[width=0.5\textwidth]{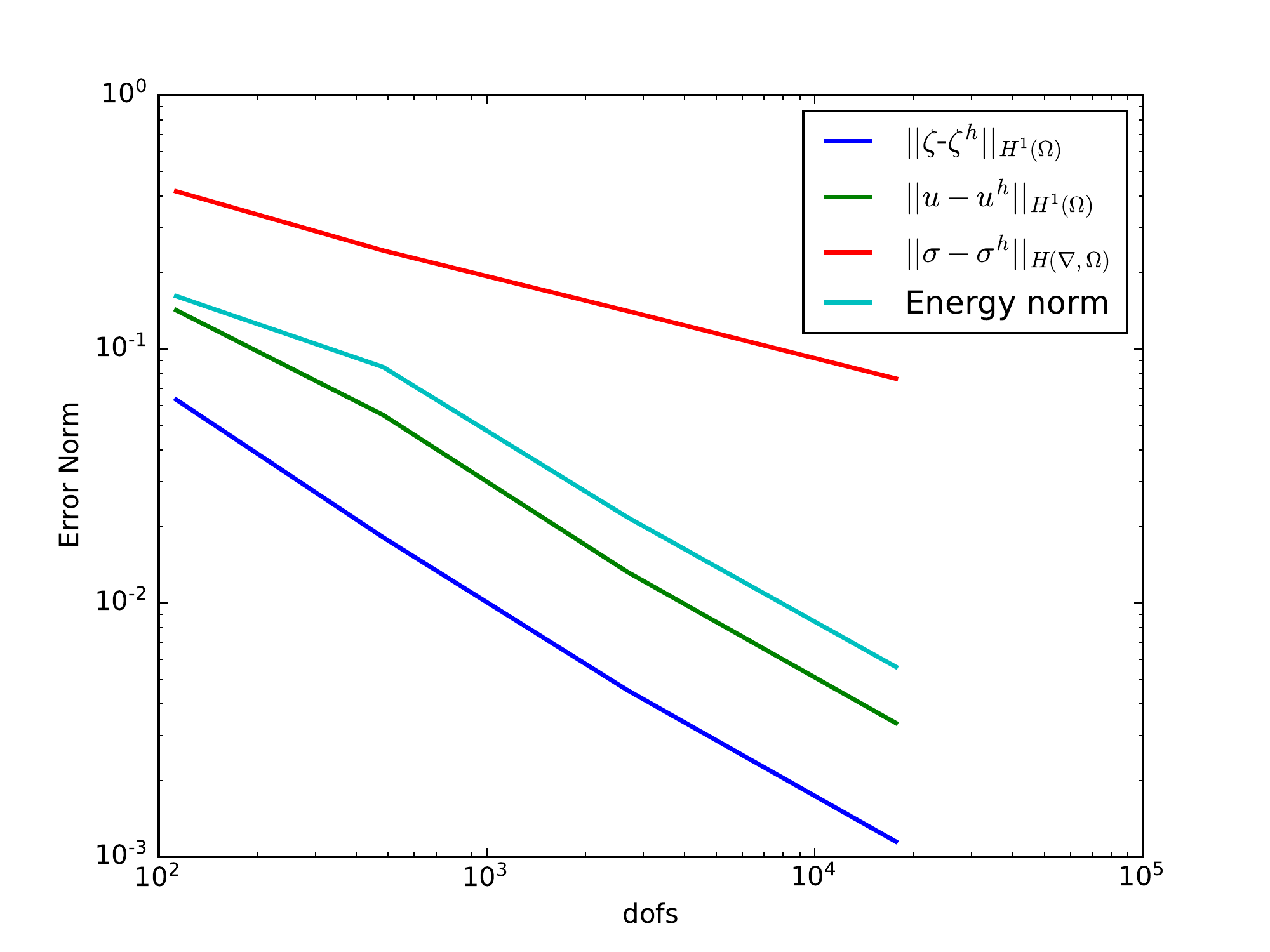}}
  \caption{\label{fig:SWE_uniform_p2} Error convergence results for uniform $h$-refinements for the 
 SWE using polynomial approximations that are quadratic for $\zeta^h$, $\uu^h$ and linear for
$\sig^h$. }
\end{figure}
In Figure~\ref{fig:SWE_uniform_p2}, we present the convergence results in several error norms, with the
observed rates of convergence listed in Table~\ref{tab:converg_rates}.
\begin{table}[h]
\centering
\caption{\label{tab:converg_rates} Error convergence rates  corresponding to 
Figure~\ref{fig:SWE_uniform_p2}. }
\begin{tabular}{@{}lc@{}}
\toprule
{Norm \hspace{6mm}} & { Observed rate \hspace{4mm}}   \\
\midrule \midrule

$\norm{\zeta - \zeta^h}{\SLTOT}$ & $3$      \\
$\norm{\uu - \uu^h}{\SLTOT}$& $3$      \\
$\norm{\sig - \sig^h}{\SLTO}$ & $2$      \\
$\norm{(\zeta-\zeta^h,\uu-\uu^h,\sig-\sig^h)}{\SLTOT}$ & $2$      \\
$\norm{\zeta - \zeta^h}{\SHOOT}$ & $2$      \\
$\norm{\uu - \uu^h}{\SHOOT}$ & $2$      \\
$\norm{\sig - \sig^h}{\SHdivO}$ & $1$      \\
$\norm{(\zeta-\zeta^h,\uu-\uu^h,\sig-\sig^h)}{\UUUT}$ & $1$      \\
$\norm{(\zeta-\zeta^h,\uu-\uu^h,\sig-\sig^h)}{\text{B}}$ & $2$      \\

\bottomrule
\end{tabular}
\end{table}

The rates of the energy norm, $\norm{(\zeta-\zeta^h,\uu-\uu^h,\sig-\sig^h)}{\SLTOT}$ and $\norm{(\zeta-\zeta^h,\uu-\uu^h,\sig-\sig^h)}{\UUU}$ are as predicted by the \emph{a priori} estimates of 
Section~\ref{sec:avs-estimates}. 
Also note that the rates of the individual error norms are of the same order as their underlying 
polynomial interpolants. In particular we point out that the rate observed for $\norm{\zeta - \zeta^h}{\SLTOT}$ is an order higher than expected from an Aubin-Nitche lift (see
Proposition~\ref{prp:rate_l2}).
We also point out that the error in the energy norm converges at a higher rate than 
 $\norm{\sig-\sig^h}{\SHdivO}$. Inspection of the proof of Lemma~\ref{lem:energy_bound_riesz}
reveals that the error in the energy norm does not the order of the Raviart-Thomas/$C^0$ approximations
used for $\sig^h$ but rather the polynomial degree of the error representation function/ optimal test function.
We have also performed verification for the time slice approach of Section~\ref{sec:time_slices} in 
which we perform uniform  mesh refinements on each slice and note that the results are essentially 
indistinguishable from those presented in Figure~\ref{fig:SWE_uniform_p2} for this exact solution.

\subsection{Adaptive Mesh Refinement}
\label{sec:adaptivity}

In this section, we present several numerical verifications for adaptive mesh refinement. 
We employ the built-in error indicators of the AVS-FE method~\eqref{eq:enrr_ind} and the 
adaptive strategy of D{\"o}rfler~\cite{dorfler1996convergent}. 

\begin{figure}[h!]
  \subfigure[ \label{fig:exaxt_ex2_ux} $u^{ex}_x(x,y)$.]{\centering
 \includegraphics[width=0.5\textwidth]{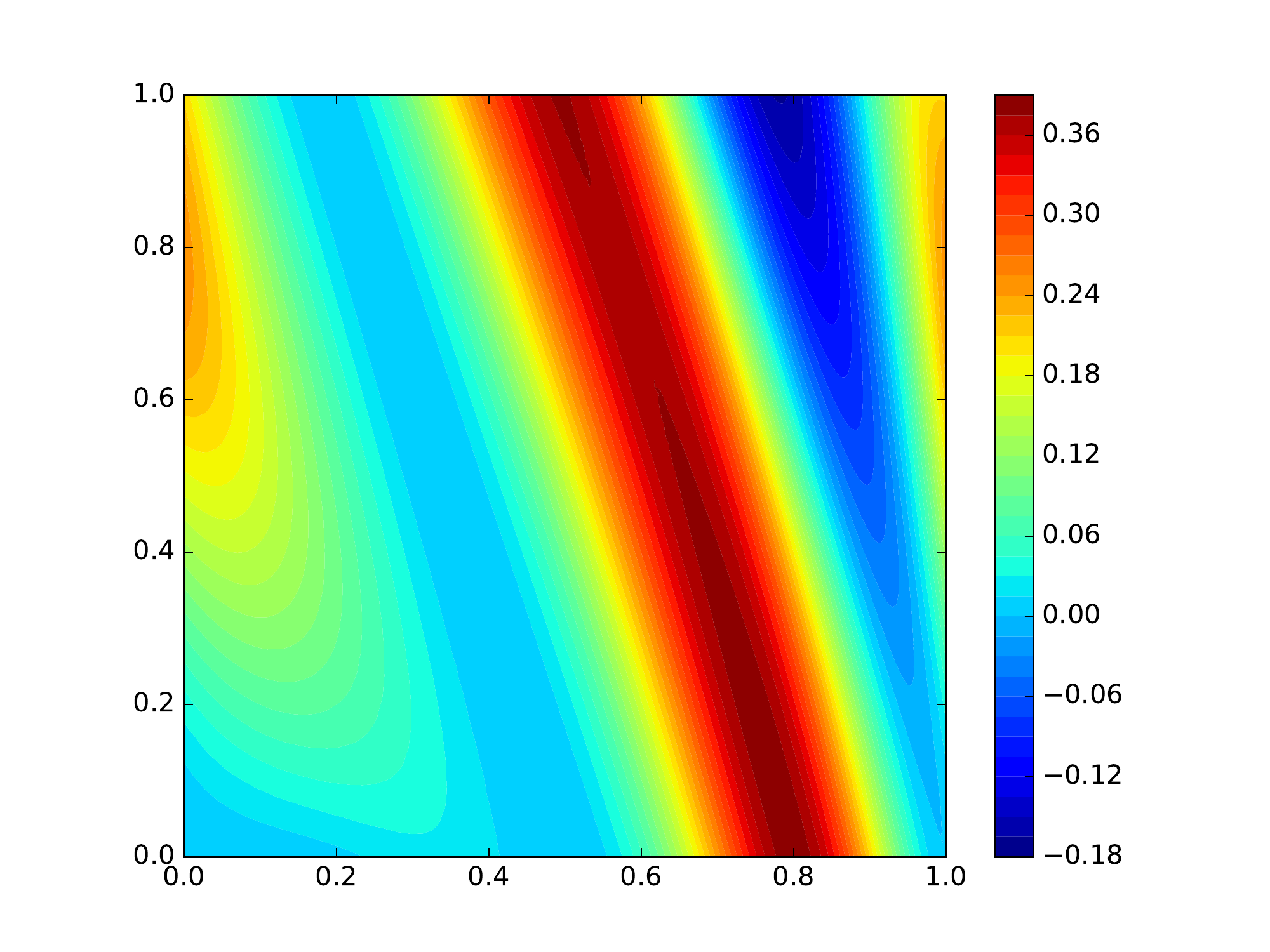}}
\hspace*{0.05in}  \subfigure[ \label{fig:exaxt_ex2_zeta} $\zeta^{ex}(x,y)$. ]{\centering
 \includegraphics[width=0.5\textwidth]{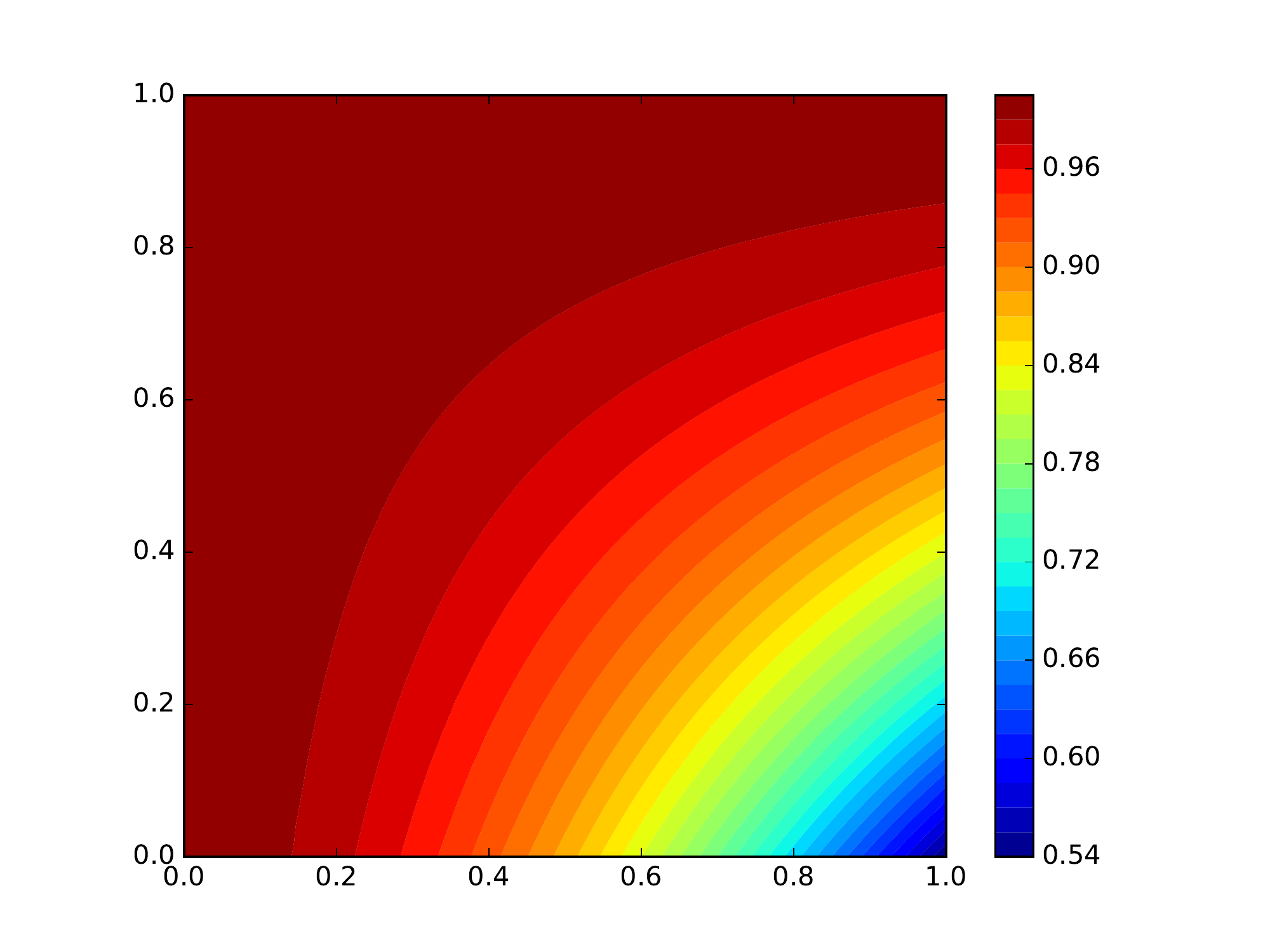}}
  \caption{\label{fig:exaxt_ex2} Exact solution components. }
\end{figure}

First, we consider a 
stationary problem on the domain $\Omega = (0,0)\times(1,1)$, and an exact solution given by:
\begin{equation} \label{eq:ex2_Exact}
\begin{array}{l}
\ds \zeta^{ex}(x,y) = \text{cos} (x - xy), \\[0.1in]
\ds u^{ex}_x(x,y) = \text{cos}^2 (\pi x + y)  \text{sin}^2 (\pi x^3 + y), \\[0.1in]
\ds u^{ex}_y(x,y) = u^{ex}_x(x,y).
 \end{array}
\end{equation}
We elect to consider a stationary problem here to ensure a meaningful presentation of the 
resulting adapted mesh.
In Figure~\ref{fig:exaxt_ex2} the components of this exact solution are presented. The 
The initial mesh consists of 2 equal triangular elements, we pick continuous quadratic polynomial
bases for $\zeta^h$, $\uu^h$, first order Raviart-Thomas bases for $\sig^h$, and
discontinuous quadratic polynomial bases for the components of the error representation function. 
The physical parameters in this case are $\mu = 10^{-9}$ and $\tau_{bf} = 1$.
\begin{figure}[h!]
  \subfigure[ \label{fig:stat_adapt_con} Convergence history.]{\centering
 \includegraphics[width=0.5\textwidth]{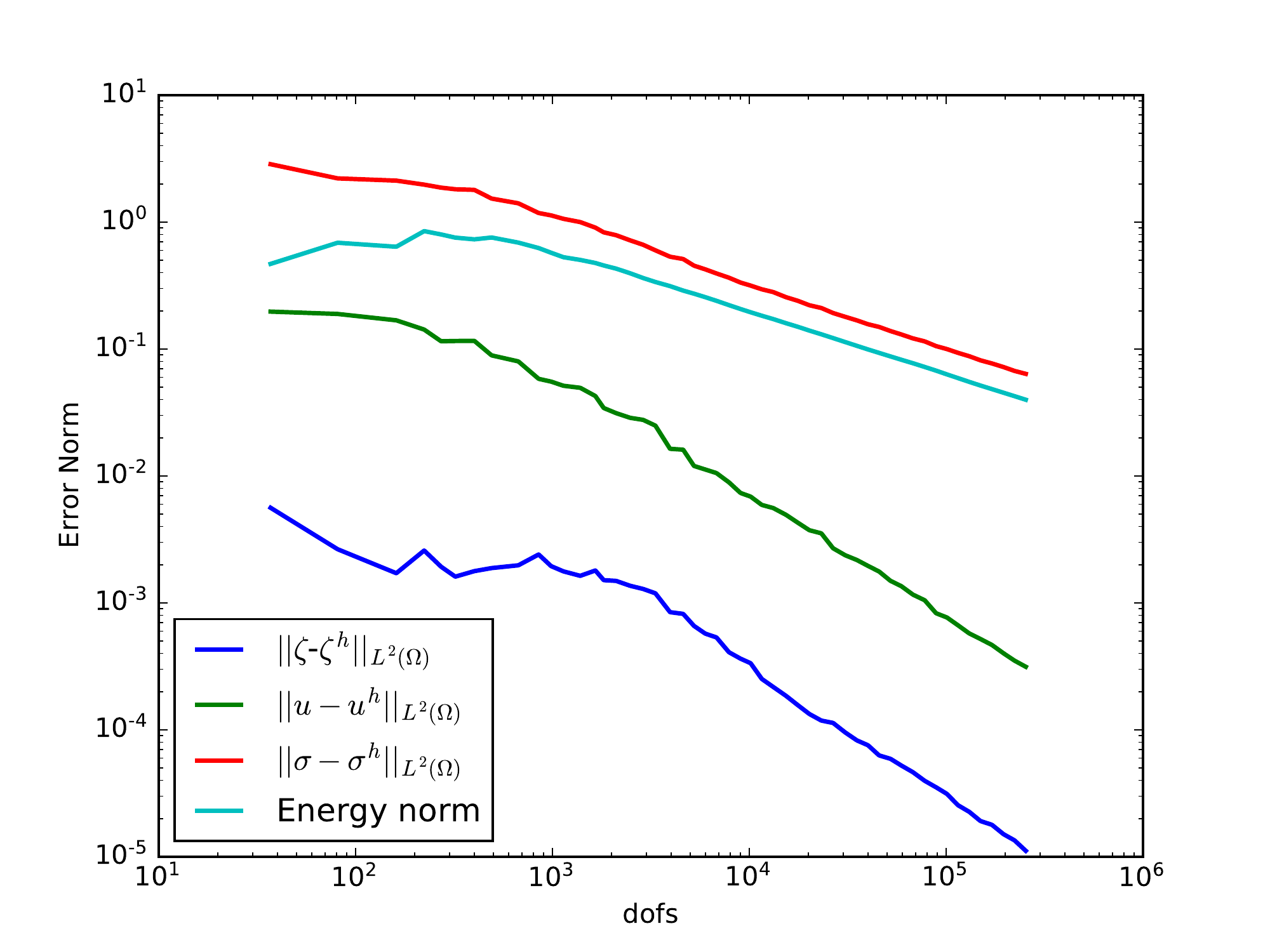}}
\hspace*{0.05in}  \subfigure[ \label{fig:adapt_fin_mesh} Final adapted mesh. ]{\centering
 \includegraphics[width=0.5\textwidth]{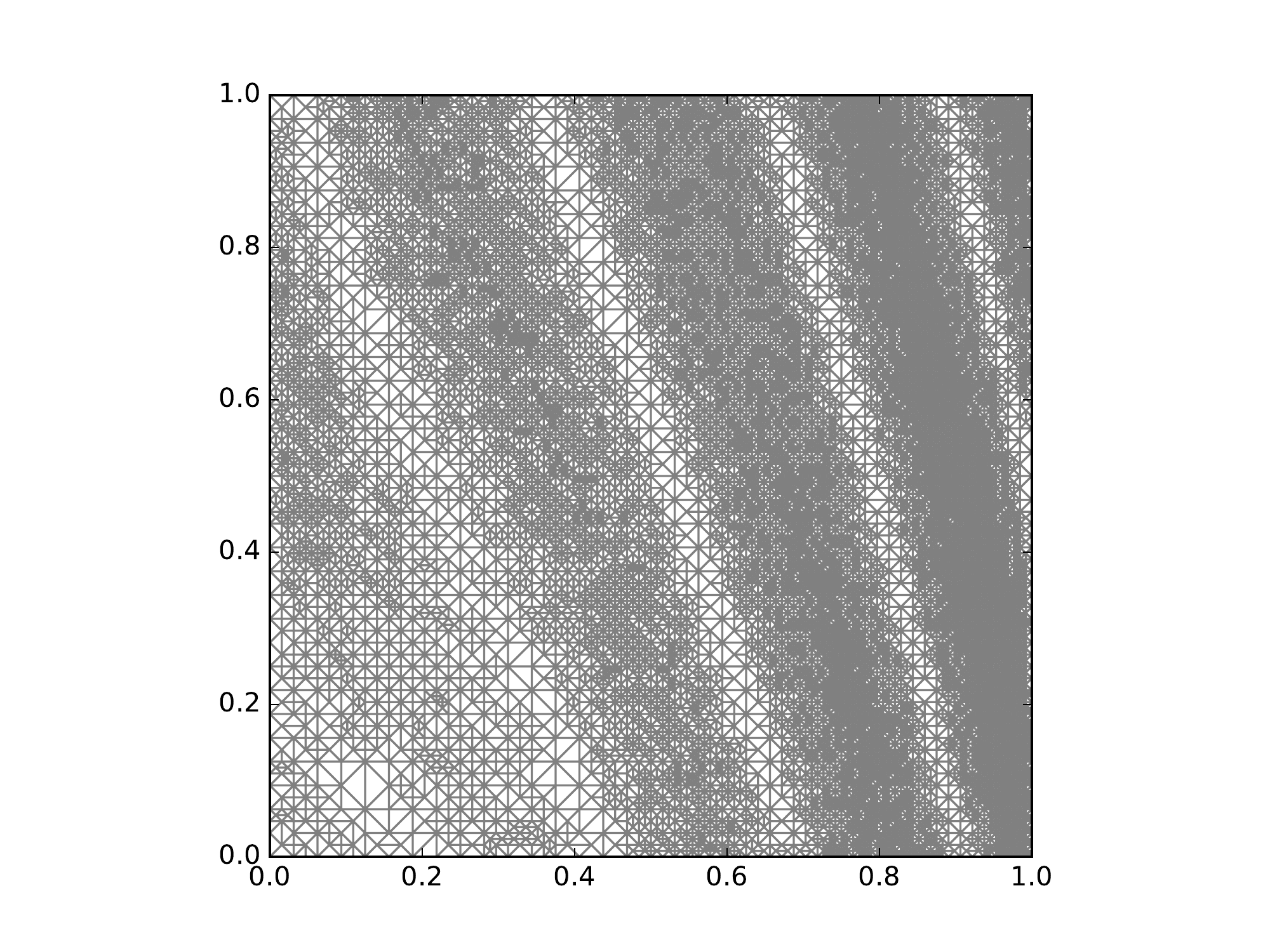}}
  \caption{\label{fig:stat_adapt} Results of the adaptive scheme applied to a stationary problem. }
\end{figure}
In Figure~\ref{fig:stat_adapt}, the final adapted mesh, after 50 refinements, is shown along with the 
convergence history of the refinement process. Comparison of mesh in Figure~\ref{fig:adapt_fin_mesh}
and the exact velocity component $u^{ex}_x(x,y)$ in Figure~\ref{fig:exaxt_ex2_ux} shows 
that the refinement process leads to a mesh that is highly refined in the regions of high velocity magnitude.

\begin{figure}[h!]
\subfigure[ \label{fig:SWE_adaptp2_l2} $\SLTOT$ and Energy error norms. ]{\centering
 \includegraphics[width=0.5\textwidth]{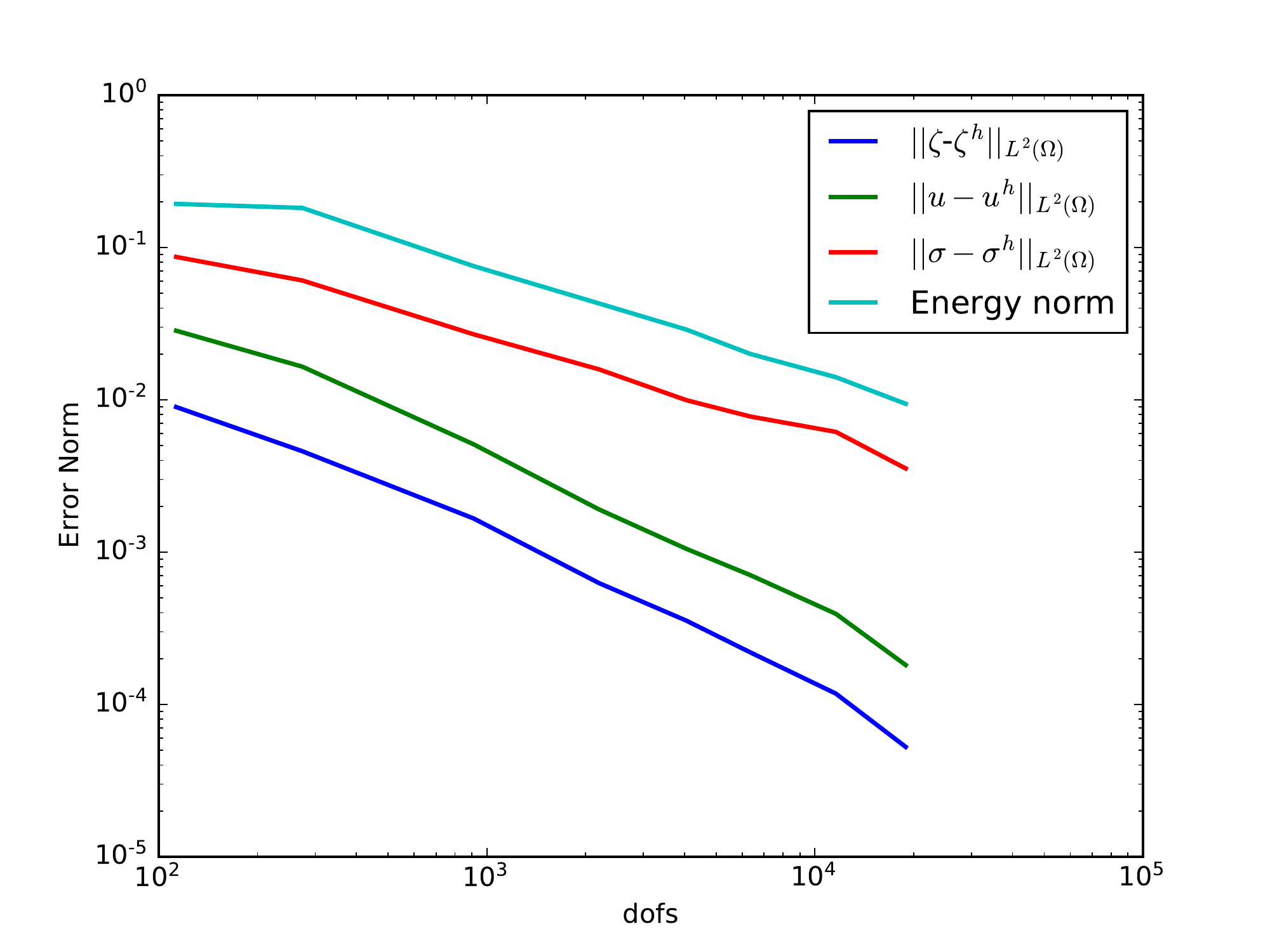}}
  \subfigure[ \label{fig:SWE_adapt_p2_h1} $\SHOOT$, $\SHdivO$ and Energy error norms.]{\centering
 \includegraphics[width=0.5\textwidth]{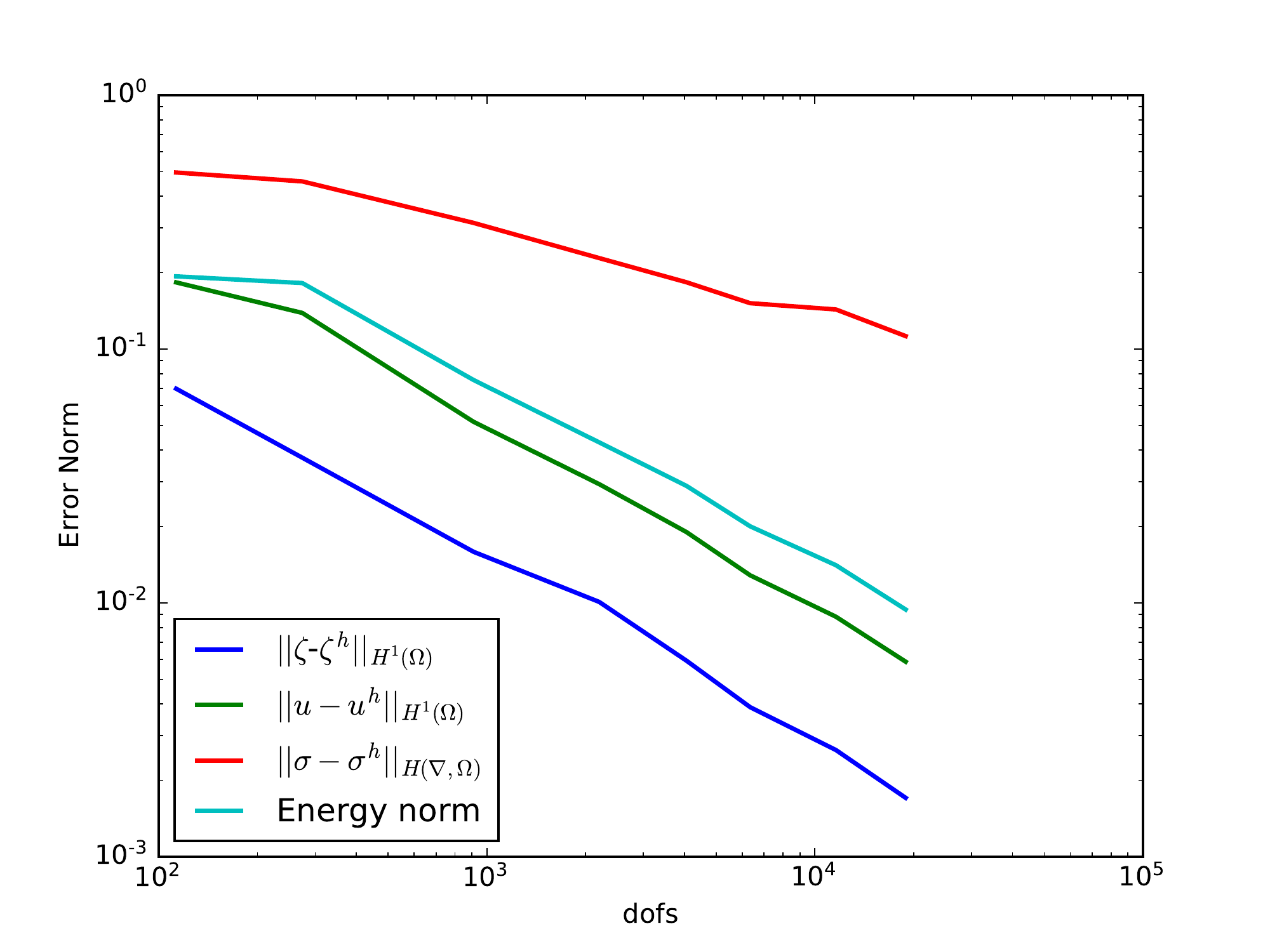}}
  \caption{\label{fig:SWE_adapt_p2} Error convergence results for adaptive $h$-refinements for the 
  purely convective flow regime. }
\end{figure}

To highlight the capabilities of the AVS-FE method of performing local time stepping, i.e., space-time
mesh adaptivity, we consider the 
limiting case for the SWE of purely convective flow. Hence, we set $\mu = 0$, $\tau_{bf} = 1$, and the
final time $T=1.0s$, and we consider the exact solutions given
in~\eqref{eq:u_exact_ex1} and~\eqref{eq:zeta_exact_ex1}. 
The initial mesh consists of six space-time tetrahedrons using quadratic polynomial basis functions 
for $\zeta^h$, $\uu^h$, linear polynomials for $\sig^h$, and
discontinuous quadratic polynomial bases for the error representation function. We perform a total of
eight mesh refinements, and report convergence histories in Figure~\ref{fig:SWE_adapt_p2}.
Space-time mesh adaptive refinements using the built-in error estimate are able to 
appropriately refine the mesh to minimize the AVS-FE numerical approximation errors as indicated 
by the results in Figure~\ref{fig:SWE_adapt_p2}. 

Next, we consider a verification of the time slice approach introduced in 
Section~\ref{sec:time_slices}, 
we again consider the problem as in the preceding verification, $\mu = 10^{-5}$, $\tau_{bf} = 1$,  
 with  a final time of $4s$.
To compare the full space-time and the space-time slices, we solve the problem using both approaches 
employing adaptive mesh refinements according to the built-in error indicators~\eqref{eq:enrr_ind} 
and the same D{\"o}rfler marking strategy. For the full space-time method, the initial mesh consists 
of 12 tetrahedron elements and  we perform a total of six mesh refinements. 
Conversely, for the time slice approach, we partition the space-time domain into eight equal slices
discretized using 12 tetrahedron elements.  In both cases, we use
quadratic polynomial basis functions for $\zeta^h$, $\uu^h$, linear polynomials for $\sig^h$, and
discontinuous quadratic polynomial bases for the error representation function. 
\begin{figure}[h!]
\subfigure[ \label{fig:SWE_slic_l2z} $\norm{\zeta-\zeta^h}{\SLTOT}$. ]{\centering
 \includegraphics[width=0.4\textwidth]{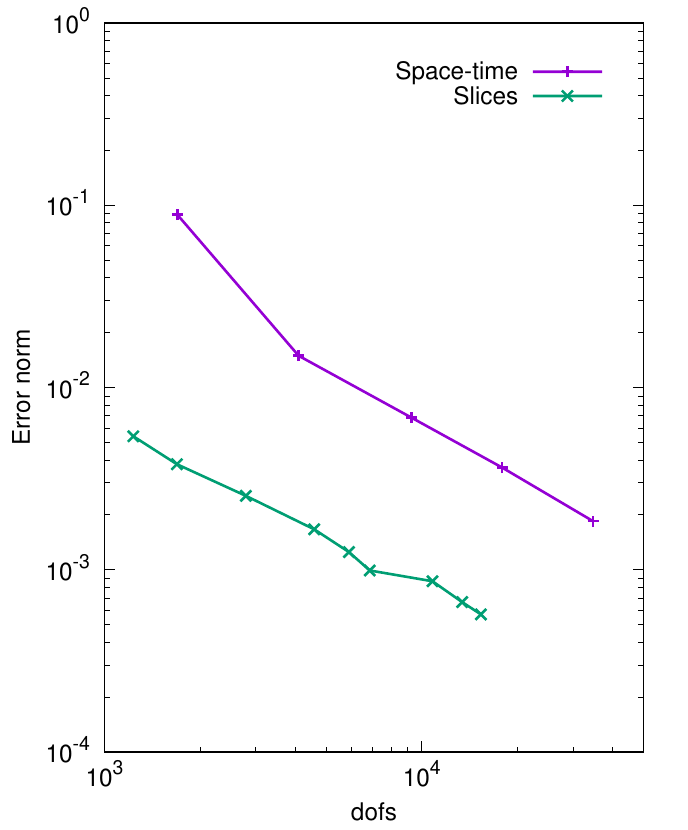}}
  \subfigure[ \label{fig:SWE_slic_l2u} $\norm{\uu-\uu^h}{\SLTOT}$.]{\centering
 \includegraphics[width=0.4\textwidth]{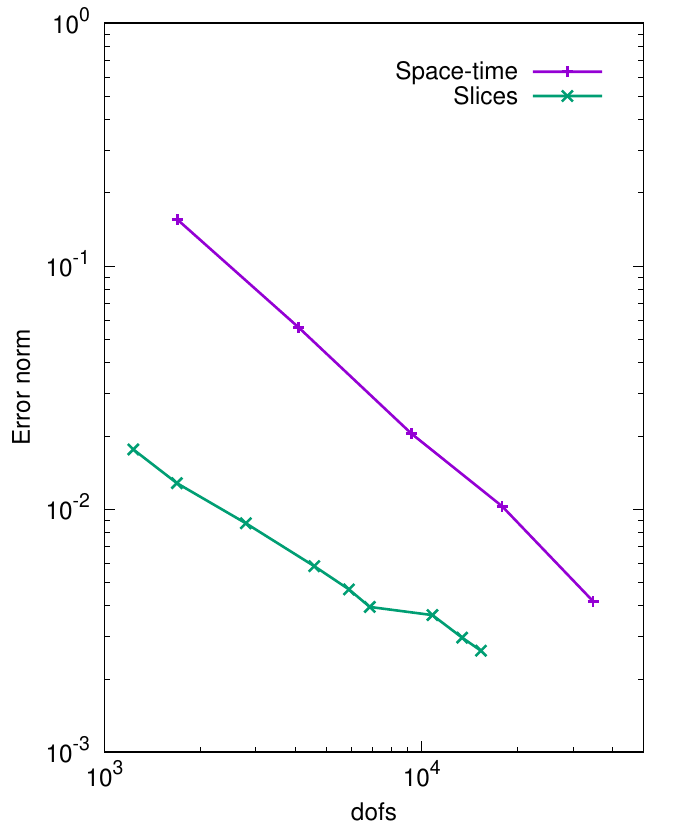}}
  \caption{\label{fig:slice_V_ST} Error convergence results comparing the space-time and time slice approaches. }
\end{figure}
\begin{figure}[h!]
\subfigure[ \label{fig:SWE_slic_l2sig} $\norm{\sig-\sig^h}{\SLTOT}$. ]{\centering
 \includegraphics[width=0.4\textwidth]{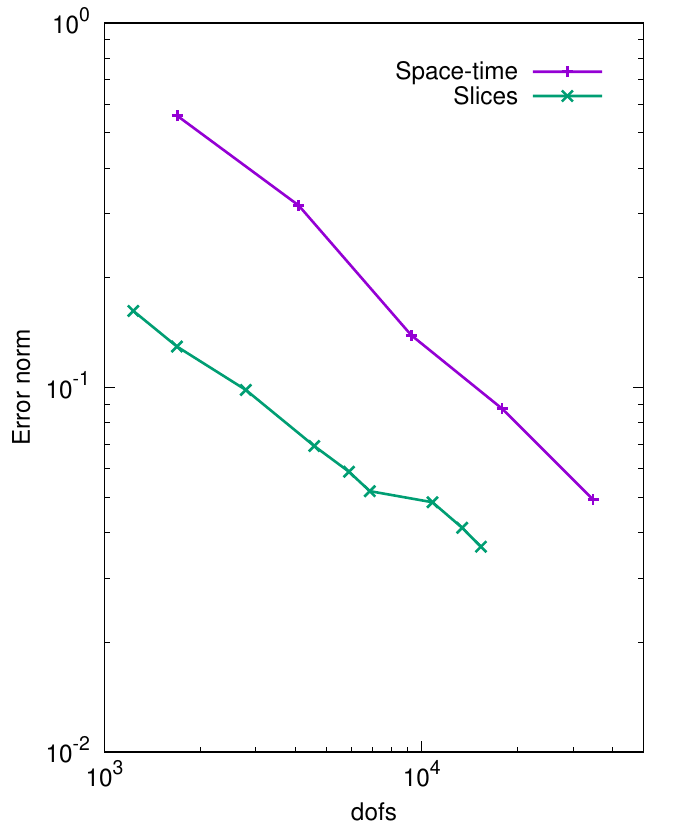}}
  \subfigure[ \label{fig:SWE_slic_E} $\norm{(\zeta-\zeta^h,\uu-\uu^h,\sig-\sig^h)}{\rm{B}}$.]{\centering
 \includegraphics[width=0.4\textwidth]{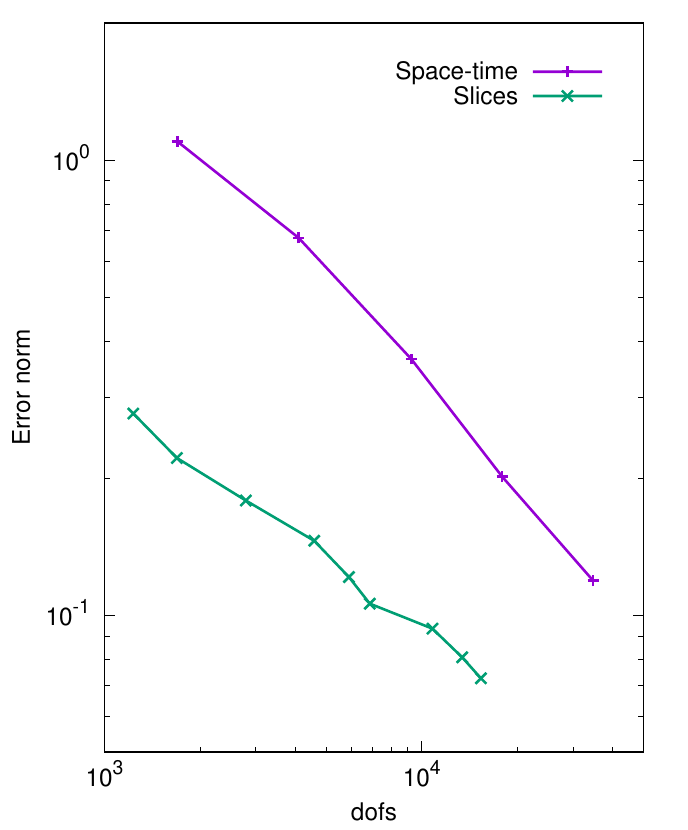}}
  \caption{\label{fig:slice_V_ST1} Error convergence results comparing the space-time and time slice approaches. }
\end{figure}
Inspection of Figure~\ref{fig:slice_V_ST1} and~\ref{fig:slice_V_ST} shows that in this case, the 
time slice approach is about an order of magnitude more accurate for the $L^2$ errors in $\zeta$ and 
$\uu$. Furthermore, the 
number of degrees of freedom is also significantly lower for the slices than for the full 
space-time approach. We have observed that for problems where the final time is less than one second, the two approaches yield essentially the same accuracy. However, as this example shows, once the final 
time becomes larger, the time slice approach is superior.

\subsection{Lake At Rest}
\label{sec:lake_at_rest}
In cases of variable bathymetry $h_b(\xx)$, a common concern in the approximation of the SWE is the 
ability to preserve the steady state of a lake-at-rest throughout the time stepping procedures 
commonly employed. Numerical approximations that preserve this steady state are referred to as 
being well balanced~\cite{michel2017well,leveque1998balancing}.
We consider a two-dimensional case $\Omega = (0,1m)\times(0,1m)$, with physical parameters 
$\mu = 10^{-5} m$ and $\tau_b = 1 s^{-1}$, and boundary conditions:
\begin{equation} \label{eq:ex4a_Set_up}
\begin{array}{l}
\ds \zeta_0 = 0, \\[0.1in]
\ds u_0 = 0, \\[0.1in]
\ds u = 0, \text{on} \, \partial \Omega, \\[0.1in]
\ds \zeta = 1m, \text{on} \, \partial \Omega. \\[0.1in]
 \end{array}
\end{equation}
Consequently, there are no sources to induce flow in this physical example. 
The bathymetry is $h_b(\xx) = 2m-h_0(\xx)$,
where $h_0$ is:
\begin{equation} \label{eq:ex_4a_bath1}
    h_0(\xx)=
    \begin{cases}
      62.5 \frac{1}{m^3}(x-0.3m)(x-0.7m)(y-0.3m)(y-0.7), & \text{if}\, (0.3m <x<0.7m)  \wedge (0.3m <y<0.7m) \\
      0, & \text{otherwise}
    \end{cases}.
\end{equation}
The final time is set to $10s$. 
We consider a space-time mesh partition consisting of 150 space-time tetrahedrons such that the "width" 
of the elements in the time direction are equal to $10s$. The resulting errors in elevation and 
velocity are shown in Table~\ref{tab:lake_At_rest}. The $L^2$ errors are vanishingly small, leading us 
to conclude that this scheme is well balanced.
\begin{table}[h]
\centering
\caption{\label{tab:lake_At_rest} Elevation and velocity errors for the lake at rest.}
\begin{tabular}{@{}ll@{}}
\toprule
{$\norm{\zeta - \zeta^h}{\SLTOT}$ \hspace{6mm}} & {$\norm{\uu - \uu^h}{\SLTOT}$ \hspace{6mm}}   \\
\midrule \midrule

$9.02 \cdot 10^{-15}$ & $4.00 \cdot 10^{-13}$      \\

\bottomrule
\end{tabular}
\end{table}

\subsection{Tidal Fluctuations}
\label{sec:tidal}

An important source impacting the flows governed by the SWE are tidal forces, as the hurricane storm
 surge can be greatly increased by tidal fluctuations. Hence, it is critical 
that the AVS-FE approximations of the SWE are able to accurately reproduce phenomena corresponding 
to  tidal fluctuations. To this end, we consider a model problem from~\cite{dawson2002discontinuous}
to facilitate comparison with existing FE methods for the SWE. 
The spatial domain is a one-dimensional channel $\Omega = (x_L,x_R) = (0,10000m)$ with a constant
bathymetry $h_b(\xx) = 10m$, the physical parameters are $\mu = 25m$ and $\tau_b = 0.01 s^{-1}$,
 and the initial and boundary conditions are:
\begin{equation} \label{eq:ex4_Set_up}
\begin{array}{l}
\ds \zeta_0 = 0, \\[0.1in]
\ds u_0 = 0, \\[0.1in]
\ds \sig(0,t) = 0, \\[0.1in]
\ds u(10000m,t) = 0, \\[0.1in]
\ds \zeta(0,t) = 0.1 \text{cos} (t \alpha) \, m,
 \end{array}
\end{equation}
where $\alpha=0.00014051891708$. Since this is a one-dimensional problem, the spaces $\SHdivO$ and
 $H^1(\Omega)$ coincide and we use $C^0$ polynomial approximations for all trial variables. 
As the period of tidal fluctuations are on the order of days, we consider a case in which the 
fluctuations occur over $7$ days, i.e., the space-time domain is 
$\Omega_{\text{T}} = (0,10000m) \times (0,604800s)$. In the corresponding AVS-FE discretization of the
space-time domain we 
employ a uniform mesh of $2(25\times400)$ triangular elements, which corresponds to a 
"time step" of 1512 seconds.
\begin{figure}[h!]
\centering
 \includegraphics[width=0.4\textwidth,angle=-90,origin=c]{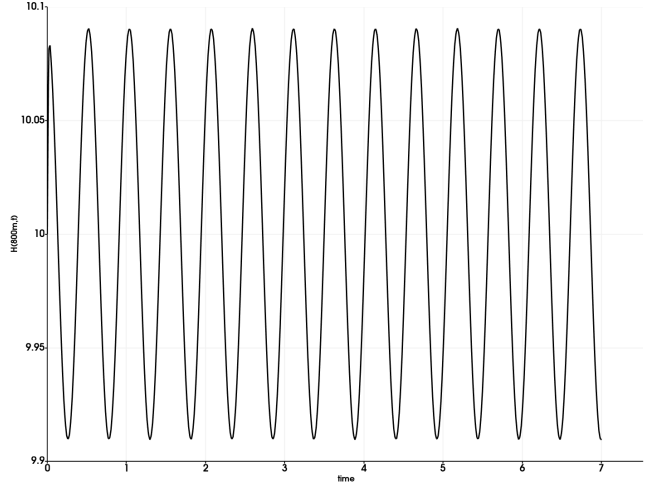}
 \vspace{-0.5in}
  \caption{\label{fig:tidal} $H^h(800m,t)$ ($m$). }
\end{figure}
In Figure~\ref{fig:tidal}, the water column elevation at $x=800m$ is shown for the full 7 day time span. As expected, the resulting 
tidal fluctuation leads to a sinusoidal elevation profile and at $x=800m$, the wave amplitude is
slightly damped from the incoming tidal wave. 
\begin{figure}[h!]
\centering
 \includegraphics[width=0.4\textwidth,angle=-90,origin=c]{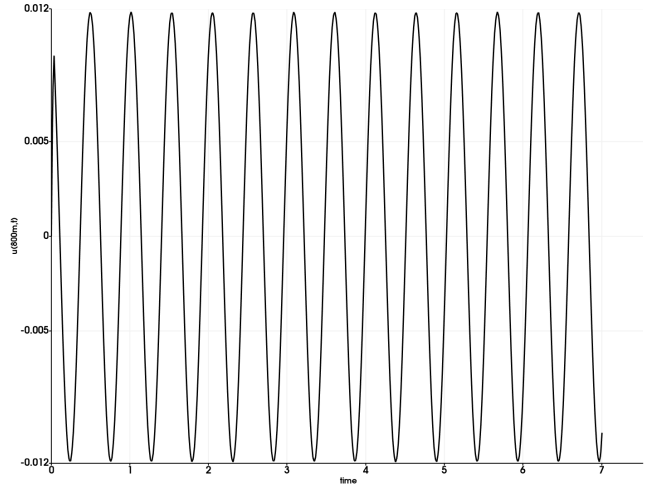}
  \vspace{-0.5in}
  \caption{\label{fig:tidalvel} $u^h(800m,t)$ ($\frac{m}{s}$). }
\end{figure}
Correspondingly, in Figure~\ref{fig:tidalvel}, the corresponding velocity is presented. Here we observe 
that the highest velocity magnitudes occur slightly before and after the peak tidal elevation.  
Visual inspection of the results in~\cite{dawson2002discontinuous}, where a coupled discontinuous and continuous Galerkin method is considered, compared with 
Figures~\ref{fig:tidal} and~\ref{fig:tidalvel} shows comparable behavior. However, the time step used 
in~\cite{dawson2002discontinuous} is $0.25s$ in Euler's method, i.e., the total number of time steps 
over seven days is roughly $2.4$ million, whereas the AVS-FE solution is obtained in a single solution 
step.

\subsection{Dam Break}
\label{sec:dam_break}

\begin{figure}[h!]
\centering
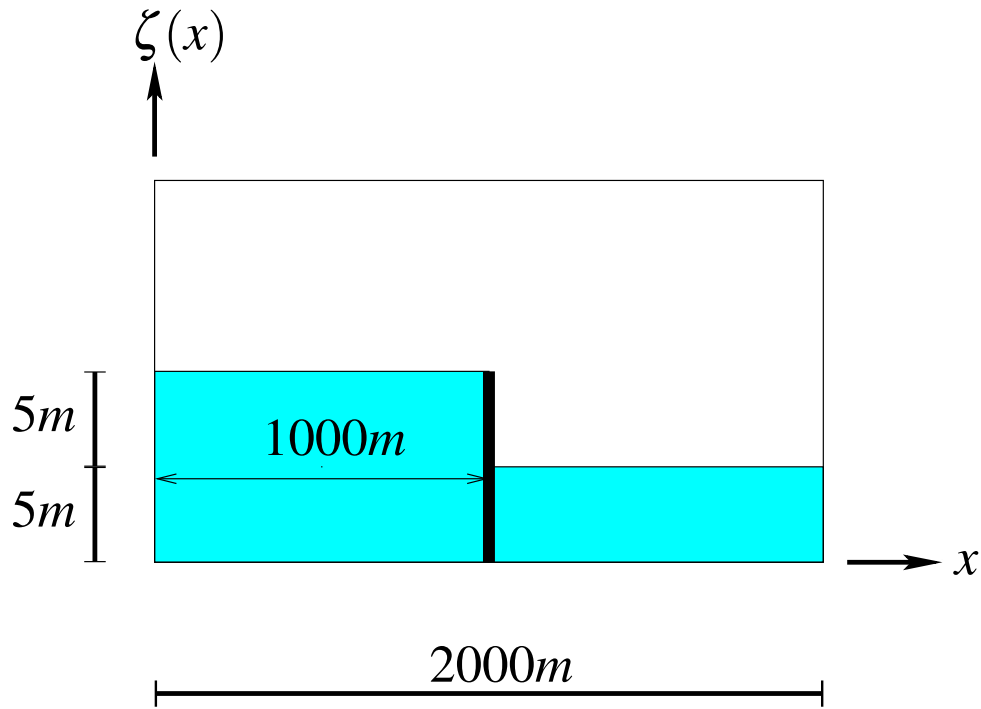 
\caption{Dam break problem spatial domain.}
\label{fig:dam}
\end{figure}

As a final numerical verification, we consider another commonly applied benchmark problem found in
literature~\cite{jacobs2015firedrake}, known as a dam break problem. We consider a one-dimensional case, in which a $2000m$ channel is divided 
by a dam separating two distinct water levels, i.e., $\Omega = (x_L,x_R) = (0,2000m)$. 
In Figure~\ref{fig:dam} the set up for the problem is shown. At $t=0s$, the dam is removed to 
simulate its total failure allowing the water to flow unconstrained into the lower reservoir. 
The physical parameters are $\nu = 10^{-2} m$ and $\tau_b = 1 s^{-1}$, we consider a  constant bathymetry 
$h_b(\xx) = 0$,
 and the initial and boundary conditions are:
\begin{equation} \label{eq:ex5_Set_up}
\begin{array}{l}
\ds \zeta_0 = 10m, \quad x \le 1000m  \\[0.1in]
\ds \zeta_0 = 5m, \quad x > 1000m  \\[0.1in]
\ds u_0 = 0, \\[0.1in]
\ds u(0,t) = 0, \\[0.1in]
\ds \zeta (2000m,t) = 0. \\[0.1in]
 \end{array}
\end{equation}
Hence, we expect sharp interfaces in the resulting velocity and elevation fields, i.e., shocks 
to develop. 
We consider the final time to be $T=200s$, i.e., the space-time domain is 
$\Omega_{\text{T}} = (0,2000m) \times (0,200s)$. The corresponding AVS-FE discretization employs a 
uniform mesh of $2(800\times35)$ triangular elements and quadratic polynomial approximations for all 
variables, except the stress which is linear. 
\begin{figure}[h!]
\centering
 \includegraphics[width=0.75\textwidth]{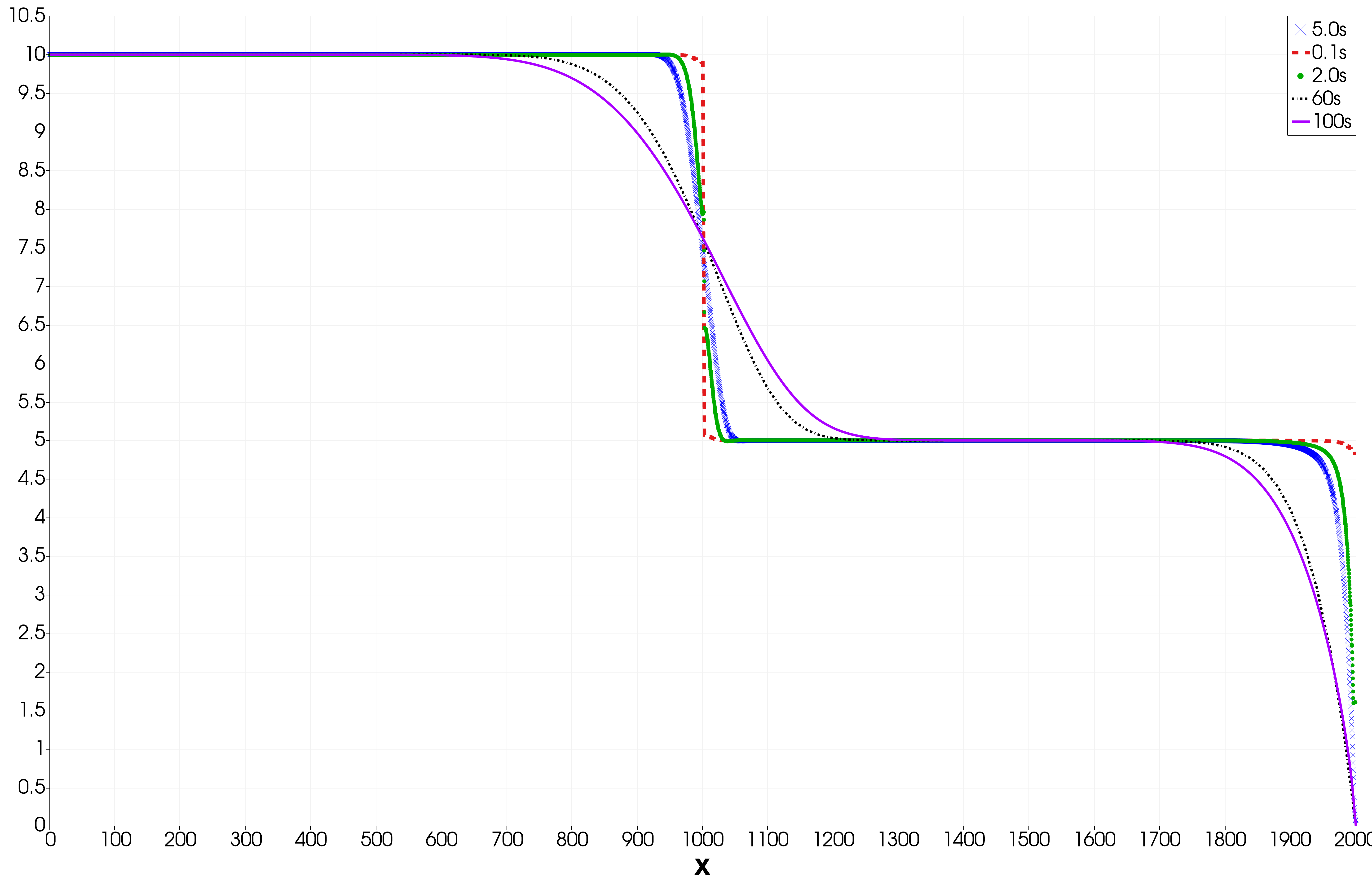}
  \caption{\label{fig:zetta} $\zeta^h(x,t)$ at various times. }
\end{figure}
In Figure~\ref{fig:zetta}, we present the elevation profile for select times. As time progresses, the
elevation profile stretches out as expected from the boundary and initial conditions. 
Shortly after the simulated dam break, at $t=0.1$ seconds, the elevation profile does not exhibit 
any noticeable oscillations leading to an accurate representation of the shock. 
Comparison with the results in~\cite{jacobs2015firedrake}, which uses slightly different boundary
and initial conditions, shows good agreement based on visual inspection.

\section{Conclusions}
\label{sec:conclusions}

We have introduced an unconditionally stable space-time FE  method for the SWE, the AVS-FE method. 
This Petrov-Galerkin method derives its stability from the DPG concept of optimal 
 test functions.  The flexibility of the DPG method allows us to establish  
continuous and stable FE approximations in both space and time by breaking the test space and 
introducing a Riesz representation problem governing the optimal test functions. 
Compared to existing FE technologies for the SWE we do not need to consider surrogate models such 
as the diffusive wave model or perform arduous problem-dependent analysis to establish 
proper stabilization parameters to achieve discrete 
stability.

Consideration of a linearized SWE allows us to establish well posedness of the AVS-FE weak formulation 
in terms of the energy norm induced by the sesquilinear form
by employing the DPG philosophy. Furthermore, for the linearized problem, we establish \emph{a priori} 
error estimates. The convergence behavior predicted by these estimates is confirmed through a sequence 
of numerical convergence studies for the full nonlinear SWE. In the case of asymptotic $h-$convergence, 
we observe optimal rates for all applicable norms of the numerical approximation errors. 
The built-in error estimate and corresponding error indicators of the AVS-FE approximation 
error in terms of the energy norm allows us to pursue space-time adaptive mesh refinement strategies.

In an effort to keep the computational cost of the space-time AVS-FE approximations low, we also 
consider a technique in which we partition the space-time domain into space-time slices. 
While the space-time AVS-FE method allows us to perform local time stepping in the form of adaptive
mesh refinements in space-time, the time slice approach allows further localization of the space-time
 mesh refinements, which now occur on each individual slice. 
We present numerical verifications in which we consider and compare the space-time and space-time slice
approaches to each other. These verifications show that the time slice approach becomes 
preferable for longer simulations. Note that what constitutes a long simulation is 
problem dependent. For the SWE and the verifications we consider here our experience indicates that the 
threshold is around one second. 

We have considered space-time adaptive mesh refinements using a built-in error estimate of the AVS-FE. 
However, we are not limited to  this type of error estimator, 
and any \emph{a posteriori} error estimation technique can be applied. In particular, 
we envision the 
use of Goal-Oriented error estimates, and their error indicators, of local quantities of interest.
 These types of estimates 
have been established in~\cite{valseth2020goal} for the AVS-FE method and linear stationary
 convection-diffusion problems. For the transient portion of the domain, goal-oriented adaptive 
 algorithms such as those developed by  Mu$\tilde{\text{n}}$oz-Matute 
\emph{et. al}~\cite{munoz2019explicit} show great potential.

The results presented in this paper serves as a proof-of-concept of the space-time AVS-FE method
applied to the SWE. 
In particular, we have presented numerical verifications of two important physical processes governed by the SWE, tidal fluctuations and that of a dam failure. 
In addition to these physical examples, we have also verified the well balanced property of the 
AVS-FE method applied to the SWE numerically in Section~\ref{sec:lake_at_rest}.
 Hence, we aim to establish a establish a new paradigm of
the application of
DPG methods in the mathematical modeling of storm surge events. 
As these physical processes are complex,  occur in complex coastal domains, and exhibit 
significant temporal variability (e.g., due to changes in wind or rainfall), the local
time stepping allowed in the space-time methods is likely to be an important factor to be 
exploited in future works. 
Coupling mechanisms such as those established by Choudhary in~\cite{choudhary2019coupled} and 
by Dawson and Proft in~\cite{dawson2002discontinuous} are likely to lead to efficient algorithms 
in the modeling of storm surge by coupling, e.g., the AVS-FE method and HDG 
methods~\cite{jones2019space,arabshahi2016space}.

\section*{Acknowledgements}
This work has been supported by the United States National Science Foundation - NSF PREEVENTS Track 2
Program, under  NSF Grant Number  1855047.


\bibliographystyle{elsarticle-num}
 \bibliography{references_eirik}
\end{document}

%% file: domain.pdf_t
\begin{picture}(0,0)%
\includegraphics{domain.pdf}%
\end{picture}%
\setlength{\unitlength}{3947sp}%
\begingroup\makeatletter\ifx\SetFigFont\undefined%
\gdef\SetFigFont#1#2#3#4#5{%
  \reset@font\fontsize{#1}{#2pt}%
  \fontfamily{#3}\fontseries{#4}\fontshape{#5}%
  \selectfont}%
\fi\endgroup%
\begin{picture}(6044,2948)(1179,-8708)
\put(3226,-6961){\makebox(0,0)[lb]{\smash{{\SetFigFont{12}{24.0}{\rmdefault}{\mddefault}{\updefault}{\color[rgb]{0,0,0}$h_b(\xx)$}%
}}}}
\put(2776,-6436){\makebox(0,0)[lb]{\smash{{\SetFigFont{12}{24.0}{\rmdefault}{\mddefault}{\updefault}{\color[rgb]{0,0,0}$\zeta(\xx)$}%
}}}}
\put(4351,-7111){\makebox(0,0)[lb]{\smash{{\SetFigFont{12}{24.0}{\rmdefault}{\mddefault}{\updefault}{\color[rgb]{0,0,0}$g$}%
}}}}
\put(1276,-7111){\makebox(0,0)[lb]{\smash{{\SetFigFont{12}{24.0}{\rmdefault}{\mddefault}{\updefault}{\color[rgb]{0,0,0}$H(\xx)=\zeta(\xx)+h_b(\xx)$}%
}}}}
\end{picture}%

%% file: spacetime.pdf_t
\begin{picture}(0,0)%
\includegraphics{spacetime.pdf}%
\end{picture}%
\setlength{\unitlength}{3947sp}%
\begingroup\makeatletter\ifx\SetFigFont\undefined%
\gdef\SetFigFont#1#2#3#4#5{%
  \reset@font\fontsize{#1}{#2pt}%
  \fontfamily{#3}\fontseries{#4}\fontshape{#5}%
  \selectfont}%
\fi\endgroup%
\begin{picture}(4237,3752)(736,-4625)
\put(2926,-3136){\makebox(0,0)[lb]{\smash{{\SetFigFont{20}{24.0}{\rmdefault}{\mddefault}{\updefault}{\color[rgb]{0,0,0}$\text{Slice 1}$}%
}}}}
\put(2926,-1936){\makebox(0,0)[lb]{\smash{{\SetFigFont{20}{24.0}{\rmdefault}{\mddefault}{\updefault}{\color[rgb]{0,0,0}$\text{Slice N}$}%
}}}}
\put(2851,-4561){\makebox(0,0)[lb]{\smash{{\SetFigFont{20}{24.0}{\rmdefault}{\mddefault}{\updefault}{\color[rgb]{0,0,0}$\Omega_0$}%
}}}}
\put(826,-2836){\makebox(0,0)[lb]{\smash{{\SetFigFont{20}{24.0}{\rmdefault}{\mddefault}{\updefault}{\color[rgb]{0,0,0}$T_{1}$}%
}}}}
\put(826,-3586){\makebox(0,0)[lb]{\smash{{\SetFigFont{20}{24.0}{\rmdefault}{\mddefault}{\updefault}{\color[rgb]{0,0,0}$t=0$}%
}}}}
\put(1567,-996){\makebox(0,0)[lb]{\smash{{\SetFigFont{20}{24.0}{\rmdefault}{\mddefault}{\updefault}{\color[rgb]{0,0,0}$t$}%
}}}}
\put(4951,-4111){\makebox(0,0)[lb]{\smash{{\SetFigFont{20}{24.0}{\rmdefault}{\mddefault}{\updefault}{\color[rgb]{0,0,0}$x$}%
}}}}
\put(751,-1711){\makebox(0,0)[lb]{\smash{{\SetFigFont{20}{24.0}{\rmdefault}{\mddefault}{\updefault}{\color[rgb]{0,0,0}$T_{final}$}%
}}}}
\end{picture}%

%% file: Dam.pdf_t
\begin{picture}(0,0)%
\includegraphics{Dam.pdf}%
\end{picture}%
\setlength{\unitlength}{3947sp}%
\begingroup\makeatletter\ifx\SetFigFont\undefined%
\gdef\SetFigFont#1#2#3#4#5{%
  \reset@font\fontsize{#1}{#2pt}%
  \fontfamily{#3}\fontseries{#4}\fontshape{#5}%
  \selectfont}%
\fi\endgroup%
\begin{picture}(5955,4313)(286,-4527)
\put(2926,-4336){\makebox(0,0)[lb]{\smash{{\SetFigFont{20}{24.0}{\rmdefault}{\mddefault}{\updefault}{\color[rgb]{0,0,0}$2000m$}%
}}}}
\put(6226,-3661){\makebox(0,0)[lb]{\smash{{\SetFigFont{20}{24.0}{\rmdefault}{\mddefault}{\updefault}{\color[rgb]{0,0,0}$x$}%
}}}}
\put(1051,-361){\makebox(0,0)[lb]{\smash{{\SetFigFont{20}{24.0}{\rmdefault}{\mddefault}{\updefault}{\color[rgb]{0,0,0}$\zeta(x)$}%
}}}}
\put(301,-3361){\makebox(0,0)[lb]{\smash{{\SetFigFont{20}{24.0}{\rmdefault}{\mddefault}{\updefault}{\color[rgb]{0,0,0}$5m$}%
}}}}
\put(301,-2761){\makebox(0,0)[lb]{\smash{{\SetFigFont{20}{24.0}{\rmdefault}{\mddefault}{\updefault}{\color[rgb]{0,0,0}$5m$}%
}}}}
\put(1876,-2911){\makebox(0,0)[lb]{\smash{{\SetFigFont{20}{24.0}{\rmdefault}{\mddefault}{\updefault}{\color[rgb]{0,0,0}$1000m$}%
}}}}
\end{picture}%